%% file: xxxsurbraid.tex
\documentclass[12pt]{article}
\usepackage{a4wide}
\usepackage{theorem}
\input{define.tex}

\title{New presentations of surface braid groups}
\author{Juan Gonz\'alez-Meneses}
\date{October, 1999}

\begin{document}
\maketitle

\begin{abstract}
  In this paper we give new presentations of the braid groups and the pure 
braid groups of a closed surface. We also give an algorithm to solve the 
word problem in these groups, using the given presentations.
\footnote{\hspace{-.65cm} Keywords: Braid - Surface - Presentation - Word 
Problem.}
\footnote{\hspace{-.65cm} {\em Mathematics Subject Classification:} 
 Primary: 20F36. Secondary: 57N05.}
\footnote{\hspace{-.65cm} Partially supported by DGES-PB97-0723 and 
by the european network TMR Sing. Eq. Diff. et Feuill.}
\end{abstract}

\section{Introduction}
Let $M$ be a closed surface, not necessarily orientable,
and let ${\cal P}=\{ P_1, \ldots, P_n \}$ be 
a set of $n$ distinct points of $M$. A {\em geometric braid} over $M$ based at 
${\cal P}$ is an $n$-tuple $\Ga=(\ga_1,\ldots, \ga_n)$ of  paths, 
$ \ga_i: \; [0,1] \longrightarrow M$, such that 
\begin{itemize}
 \item[(1)] $\ga_i(0)=P_i$ for all $i=1,\ldots, n$,

 \item[(2)] $\ga_i(1)\in {\cal P}$ for all $i=1,\ldots, n$,

 \item[(3)] $\{\ga_1(t),\ldots, \ga_n(t)\}$ are $n$ distinct points of $M$ for all 
        $t\in [0,1]$.
\end{itemize}

 For all $i=1,\ldots,n$, we will call $\ga_i$ the {\em $i$-th string} of $\Ga$.

\vspace{.3cm}
Two geometric braids based at ${\cal P}$ are said to be {\it equivalent} if there 
exists a homotopy which deforms one of them into the
other, provided that at any time we always have a geometric 
braid based at ${\cal P}$. 
We can naturally define the product of two braids as induced by 
the usual product of paths: for every $i=1,\ldots,n,\;$
we compose the string of the first braid which 
ends at $P_i$, with the $i$-th string of the second braid. 
This product is clearly well defined, and it endows the set of 
equivalence classes of braids with a group structure. 
This group is called the {\em braid group on $n$ strings over 
$M$ based at ${\cal P}$}, and is denoted by $B_n(M,{\cal P})$. 
This group does not depend, up to isomorphism, on the choice of ${\cal P}$, 
but only on the number of strings, so we may write $B_n(M)$ instead 
of $B_n(M,{\cal P})$.

\vspace{.3cm}
  A braid $\Ga=(\ga_1,\ldots, \ga_n)$ is said to be {\it pure} if 
$\ga_i(1)=P_i$ for all $i=1,\ldots,n$, that is, if all its strings are 
loops. The set of equivalence classes of pure braids forms a subgroup of 
$B_n(M,{\cal P})$ called {\em pure braid group on $n$ strings over $M$ 
based at ${\cal P}$}, and denoted $PB_n(M,{\cal P})$. Again, we may write 
$PB_n(M)$ since it does not depend on the choice of ${\cal P}$. Note that 
if $n=1$, then $B_1(M)=PB_1(M)=\pi_1(M)$, the fundamental group of $M$.

\vspace{.3cm}
  There exists an interpretation of braid groups as fundamental groups of 
some spaces, called {\em configuration spaces}. Let $F_nM$ denote the space of 
$n$-tuples of distinct points of $M$, that is, 
$F_nM= M^n \backslash \De$, where
$$
   \De = \{ (x_1,\ldots,x_n)\in M^n \; / \; x_i= x_j \mbox{ for some } 
    i\neq j \}.
$$
It is clear that $PB_n(M)\simeq \pi_1(F_nM)$. Now consider the symmetric 
group on $n$ elements, $\Si_n$. This group acts naturally
on $F_nM$ by permuting coordinates, so we can consider the 
{\em configuration space}: 
$$
    \widehat{F}_nM=F_nM/\Si_n,
$$
which can be seen as the space of embeddings of $n$ points in $M$. We 
clearly have $B_n(M)\simeq \pi_1(\widehat{F}_nM)$. 

\vspace{.3cm}
  This way to look at braids provides some useful exact sequences, derived
from fibrations. The first one comes from the covering
space map
$$
   F_nM \longrightarrow \widehat{F}_nM,
$$
with fiber $\Si_n$. It induces the following exact sequence:
\begin{equation}\label{es1}
    1\longrightarrow PB_n(M) \stackrel{e}{\longrightarrow}
       B_n(M) \stackrel{f}{\longrightarrow} \Si_n \longrightarrow 1.
\end{equation}
 The homomorphism $e$ is the natural inclusion, and $f$ maps a given braid 
to the permutation that it induces on ${\cal P}$.

\vspace{.3cm}
  Now we consider the Fadell-Neuwirth fibration (\cite{fadellneuw}):
given $1\leq m < n$, the map
$$
\begin{array}{rrcl}
  p: &            F_nM  & \longrightarrow & F_mM \\
     & (x_1,\ldots,x_n) &   \longmapsto   & (x_{n-m+1},\ldots,x_n)  
\end{array}
$$
is a locally trivial fibration with fiber
$F_{n-m}(M\backslash \{Q_1,\ldots,Q_m\})$, for any choice of the points 
$\{Q_1,\ldots,Q_m\}$. Set ${\cal P}'=\{P_2,\ldots,P_n\}$, take $m=n-1$, 
and consider $M$ different from the sphere and from the projective plane 
(so $\pi_2(M)=1$). By the long exact sequence of homotopy groups of this 
fibration, we obtain
\begin{equation}\label{es2}
    1\longrightarrow \pi_1(M\backslash {\cal P}', P_1) 
      \stackrel{u}{\longrightarrow} PB_n(M,{\cal P}) 
    \stackrel{v}{\longrightarrow} PB_{n-1}(M,{\cal P}') \longrightarrow 1.
\end{equation}
 If $\ga\in \pi_1(M\backslash {\cal P}', P_1)$, then  $u(\ga)=(\ga,e_{P_2},
\ldots,e_{P_n})$, where $e_{P_i}$ denotes the constant path on $P_i$, 
and, for $\Ga=(\ga_1,\ldots,\ga_n)\in PB_n(M,{\cal P})$, one has 
$v(\Ga)=(\ga_2,\ldots,\ga_n)$.

\vspace{.3cm}
The goal of this paper is to determine new presentations of the braid groups
of closed surfaces different from the sphere and from the 
projective plane. These presentations are much simpler than those which were
known before (\cite{scott}). Moreover, the generators and the
relations have an easy geometric interpretation. We also show that these
presentations furnish an algorithm to solve the word problem for
surface braid groups. Notice that similar presentations of the braid groups of
the sphere and of the projective plane can be found in \cite{fadellvanb}
and in \cite{vanb}, respectively.

\vspace{.3cm}
  Our work is organized as follows. In Section~\ref{statements} we state 
the results, introducing the generators and relations of our new 
presentations. Then we explain in Section~\ref{method} the method followed 
in the proofs, which we apply throughout Sections~\ref{secori} and 
\ref{secnonori}, for orientable and non-orientable surfaces, respectively.
Finally, we describe in Section~\ref{word} an algorithm to solve the 
word problem in surface braid groups.

\vspace{.3cm}
  I would like to thank Luis Paris for giving me the idea of applying 
Lemma~\ref{baselem} to surface braid groups, and also for its valuable
help in the writing of this paper.

\section{Statements}\label{statements}
The aim of this section is to state our presentations
of surface braid groups, defining the generators and
showing that the proposed relations are satisfied. We start
with the case of an oriented surface different from the sphere.

\begin{figure}[ht]
\centerline{\input{polygon.pstex_t}}
\caption{The polygon $L$ representing $M$.}
\label{polygon}
\end {figure}
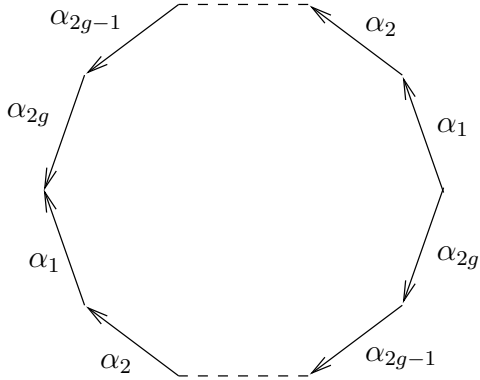

\vspace{.3cm}
 Let $M$ be a closed, orientable surface of genus $g\geq 1$. The first
thing we want is to have a geometrical representation of a braid
over $M$. We represent $M$ as a polygon $L$ of
$4g$ sides, identified in the way of Figure~\ref{polygon} 
(See \cite{massey}, page 34, ex. 8.9).

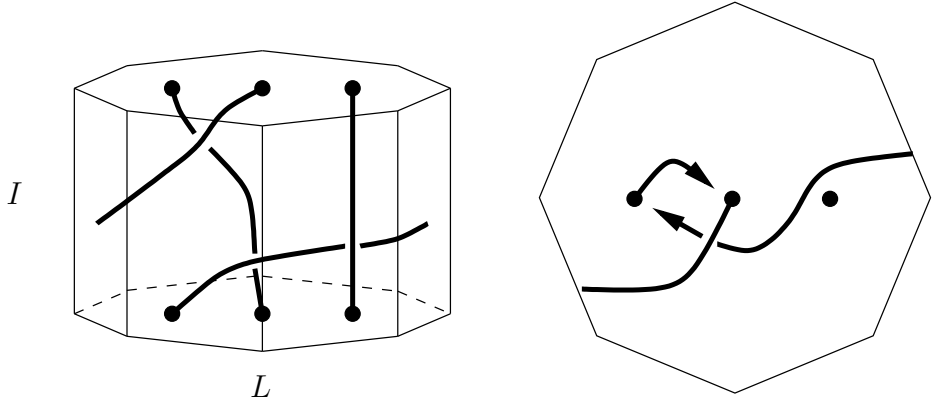
\begin{figure}[ht]
\centerline{\input{viewpoint1.pstex_t}\quad \quad \quad 
  \input{viewpoint2.pstex_t}}
\caption{A braid over a surface of genus $2$: two different viewpoints.}
\label{viewpoint}
\end {figure}

\vspace{.3cm}
  We could now take the cylinder $L\times I$ ($I=[0,1]$), and represent a 
braid $\Ga$ over $M$ as it is usually done for the open disc, that is, in 
$L\times \{t\}$ we draw the $n$ points $\ga_1(t),\ldots,\ga_n(t)$.
But in this case a string could ``go through a wall'' of the cylinder and appear 
from the other side. Hence, if we look at the cylinder from the usual viewpoint,
it would not be clear which are the ``crossed walls'' (see 
the left hand side of Figure~\ref{viewpoint}).

\vspace{.3cm}
  The solution we propose is to look at the cylinder from above, as in the 
right hand side of Figure~\ref{viewpoint}. In this way, we get rid of the 
ambiguity, and moreover we see the strings again as paths in the surface. 
When two strings cross, we see passing above the one that reaches before 
the crossing point. Anyway, it is good to keep in mind the idea that we are 
looking to a cylinder, and to consider the paths as strings: in this manner,
one can see more easily when two geometric braids are equivalent.

\begin{figure}[ht]
\centerline{\input{a1r.pstex_t}\quad \input{a1s.pstex_t}\quad 
  \input{sii.pstex_t}}
\caption{The generators of $B_n{M}$.}
\label{generators1}
\end {figure}
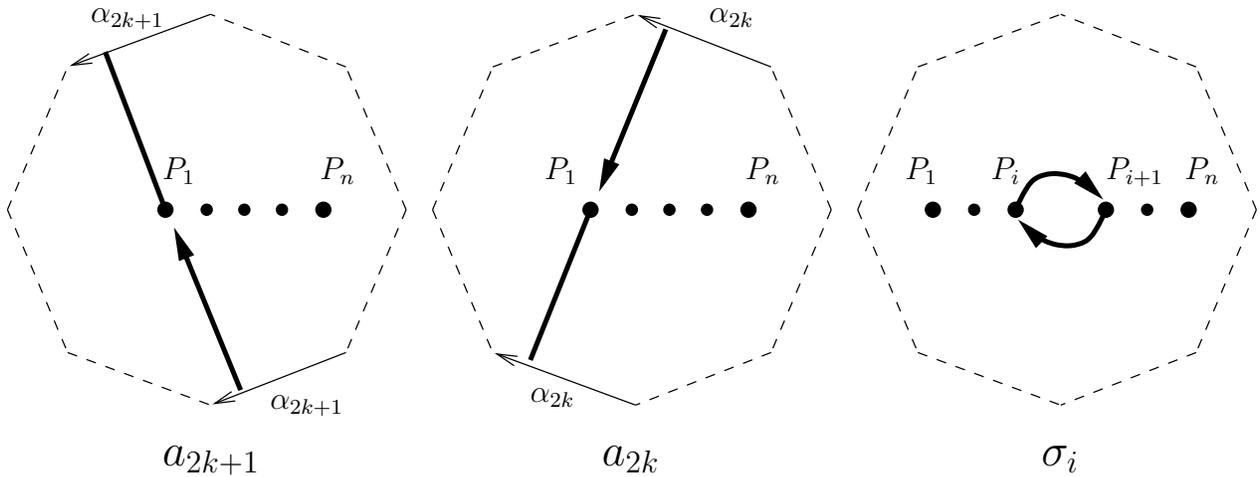

\vspace{.3cm} 
  Now we can define the generators of $B_n(M)$.
We choose the $n$ base points along the horizontal diameter of 
$L$, as in Figure~\ref{generators1}. Now given $r$, $\; 1\leq r \leq 2g$  
we define the 
braid  $a_{r}$ as follows: its only nontrivial string is
the first one, which goes through the $r$-th wall, in the way of
Figure~\ref{generators1}. That is, the first string will go upwards if $r$ is 
odd, and downwards otherwise.

\vspace{.3cm}
 We also define, for all $i=1,\ldots, n-1$, the braid $\si_i$ as in 
Figure~\ref{generators1}. Note that $\si_1,\ldots,\si_{n-1}$ are the 
classical generators of the braid group $B_n$ of the disc.

\vspace{.3cm}
  We will see later that the set $\{a_1,\ldots,a_{2g},\si_1,\ldots,\si_{n-1} \}$
is a set of generators of $B_n(M)$. There are two relations between these
generators that we can deduce as follows. Consider the interior of $L$. 
It is a subsurface $D$ of $M$ homeomorphic to a disc, so clearly every relation 
satisfied in the braid group $B_n=B_n(D)$ will be satisfied as well in $B_n(M)$ 
(the same homotopy can be used in both cases). In fact, since $g\geq 1$, it is 
known that $B_n$ is a subgroup of $B_n(M)$ (see \cite{parisrolfsen}).  
Hence, from the classical presentation of $B_n$, we obtain two relations
in $B_n(M)$: 
$$
  \si_i \si_j = \si_j \si_i  \quad  \quad \quad (|i-j|\geq 2),
$$
$$
  \si_i \si_{i+1} \si_i = \si_{i+1} \si_i \si_{i+1}
     \quad \quad \quad     (1 \leq i \leq n-2). 
$$

\vspace{.3cm}
  Note also that if $i\in \{2,\ldots,n-1\}$ and $r\in\{1,\ldots,2g\}$, then the 
non-trivial strings of $\si_i$ and the one of $a_r$ may be taken to be disjoint. 
This clearly implies that these two braids commute. Hence we have
$$
   a_{r} \si_i=\si_i a_{r}    \quad \quad \quad 
    (1\leq r \leq 2g; \; \; i\geq 2).
$$
 
\begin{figure}[ht]
\centerline{\input{puzzleP11.pstex_t} \hspace{2cm} 
        \input{puzzleP12.pstex_t}}
\caption{The initial and the $P_1$-polygon of a surface of genus $2$.}
\label{puzzleP1}
\end {figure}
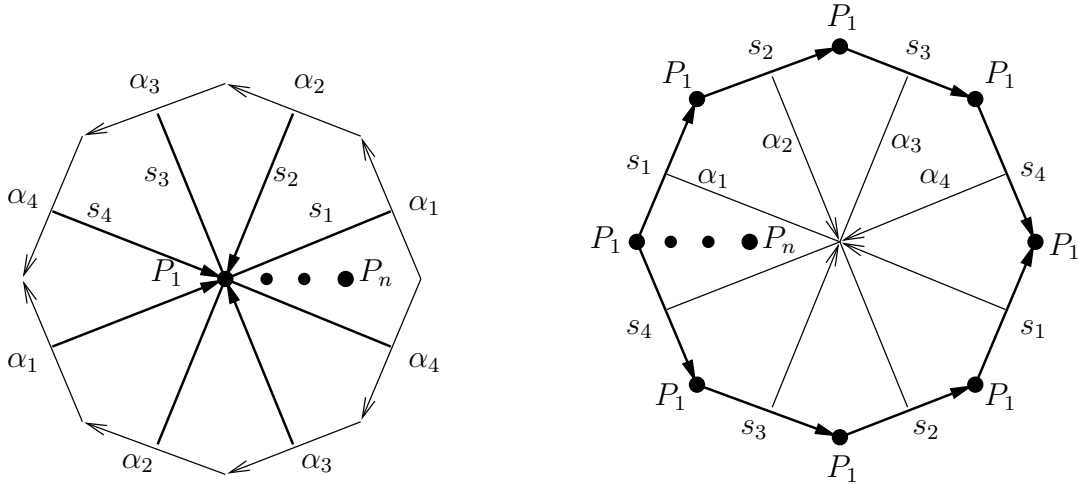

\vspace{.3cm}
  Now, in order to find more relations between the set of generators, we
do the following construction. Denote by $s_{r}$ the first string of $a_{r}$, 
for all $r=1,\ldots,2g$, and consider all the paths  $s_{1},\ldots,s_{2g}$. 
We can ``cut'' the polygon $L$ along them, and ``glue'' the pieces 
along the paths $\al_1,\ldots,\al_{2g}$. We obtain another polygon of $4g$ 
sides which are labeled by $s_{1},\ldots,s_{2g}$
(see in Figure~\ref{puzzleP1} the case of  a surface of genus $2$; the 
general case is analogous). We will call this new polygon the 
{\em $P_1$-polygon} of $M$, since all of its vertices are identified to
$P_1$, while $L$ will be called the
{\em initial polygon}. We obtain in this way a new representation of the 
surface $M$.  

\vspace{.3cm}
We will use the $P_1$-polygon to show three more relations in $B_n(M)$.
For instance, consider the braid $a_1\cdots a_{2g}a_1^{-1}\cdots a_{2g}^{-1}$.
If we look at it in the $P_1$-polygon, it is clear that it
is equivalent to the braid of Figure~\ref{B3}. But this one can be seen into
the initial polygon as a braid that does not go through the walls, namely,
an element of $B_n$, the braid group of the disc. Then
we can easily show that it is equivalent to
the braid $\si_{1}\cdots \si_{n-2} \si_{n-1}^2 \si_{n-2} \cdots \si_{1} $. 
So we have:
$$
   a_1\cdots a_{2g}a_1^{-1}\cdots a_{2g}^{-1}=
  \si_{1}\cdots \si_{n-2} \si_{n-1}^2 \si_{n-2} \cdots \si_{1}.
$$

\begin{figure}[ht]
\vspace{.3cm}
\centerline{\input{B3.pstex_t}}
\caption{The braid $a_1\cdots a_{2g}a_1^{-1}\cdots a_{2g}^{-1}$.}
\label{B3}
\end {figure}

\begin{figure}[ht]
\centerline{\input{A2r2.pstex_t}\hspace{1.5cm}\input{A2r1.pstex_t} }
\caption{The braid $A_{2,r}$: In the $P_1$-polygon and in the initial one.}
\label{A2r}
\end {figure}
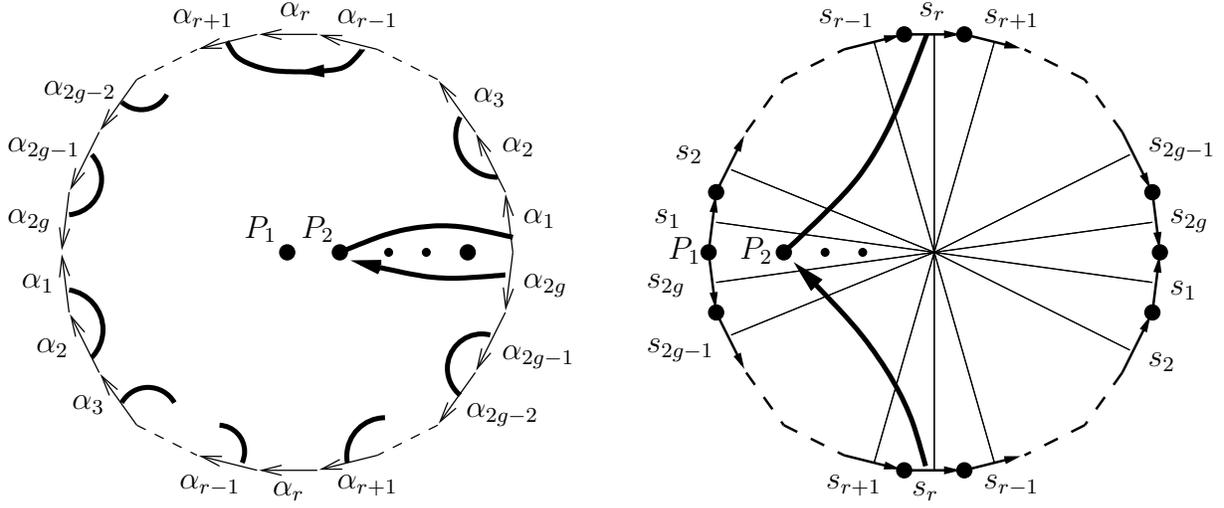

\vspace{.3cm}
  Now we define, for each $r=1,\ldots,2g$, the braid 
$$
    A_{2,r}=\si_1^{-1}
   \left( a_{1}\cdots a_{r-1}a_{r+1}^{-1}\cdots a_{2g}^{-1}\right)
    \si_1^{-1}.
$$
We will use the $P_1$-polygon to see how it looks like. 
In the left hand side of Figure~\ref{A2r}, we can see a braid which is
clearly equivalent to $A_{2,r}$ (if $r$ is odd, the other case 
being analogous). If we ``cut'' 
and ``glue'' to see this braid in the $P_1$-polygon, we obtain 
the situation of the right hand side of 
Figure~\ref{A2r}. That is, $A_{2,r}$ can be seen as a braid whose only 
nontrivial string is the second one, which goes upwards and crosses once 
the $r$-th wall $s_r$. Note that, unlike the case of $a_r$, 
$A_{2,r}$ always points upwards in the $P_1$-polygon, no matter the 
parity of $r$.

\vspace{.3cm}
Therefore we have seen that the braid $A_{2,r}$ can be 
represented by a geometric braid, whose only non trivial string 
can be taken disjoint from all the paths 
$s_t$, $\; t\neq r$. This clearly implies that
$$
    a_{t} A_{2,r}= A_{2,r} a_{t}  \quad \quad \quad
   (1\leq t,r \leq 2g; \;\;  t \neq r ).  
$$

\begin{figure}[ht]
\centerline{\input{B51.pstex_t} \quad \input{B52.pstex_t}\quad 
   \input{B53.pstex_t} }
\caption{The braid $[a_1\cdots a_{r},A_{2,r}]$.}
\label{B5}
\end {figure}
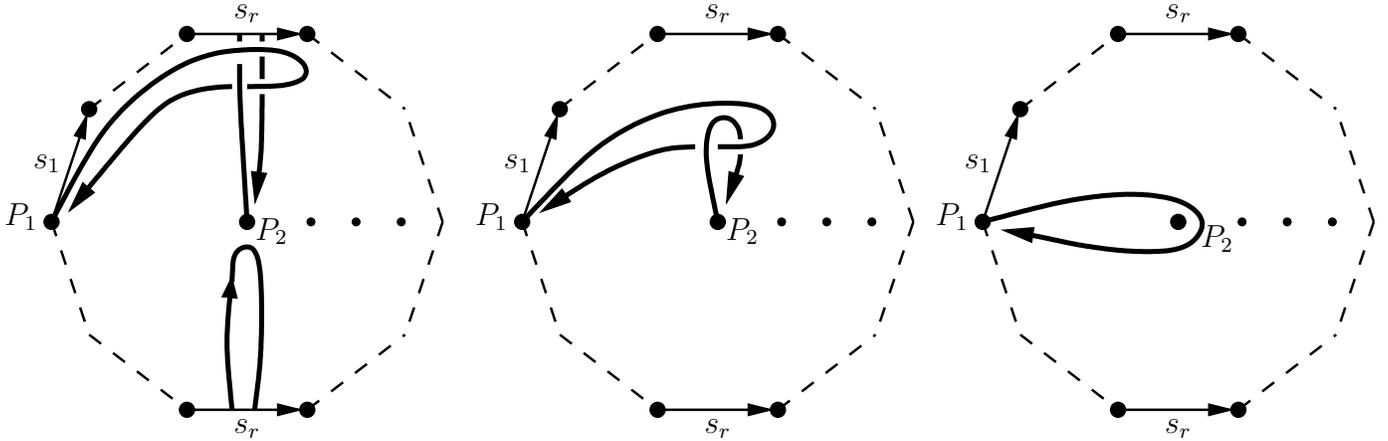

\vspace{.3cm}
  Now we finish our set of relations by considering the commutator of
the braids $\left( a_1\cdots a_{r}\right)$ and $A_{2,r}$, for all 
$r=1,\ldots,2g$.
In Figure~\ref{B5} we can see a sketch of the homotopy which starts with 
this commutator and deforms it to a braid clearly equivalent to $\si_1^2$.
Therefore, we obtain the relation:
$$
      \left( a_{1}\cdots a_{r}\right) A_{2,r}=\si_1^2 A_{2,r}
      \left( a_{1}\cdots a_{r}\right)  \quad \quad \quad 
       (1\leq r\leq 2g).     
$$

\vspace{.3cm}
  Now we claim that the six relations that we have considered form a 
complete set of defining relations of $B_n(M)$. 
In other words, we have the following result.

\vspace{.3cm}
\begin{theorem}\label{presBnM}
  If $M$ is a closed, orientable surface of genus $g\geq 1$, then 
$B_n(M)$ admits the following presentation:

\begin{itemize}

 \item Generators:
  \begin{itemize}

  \item[]   $\si_1, \ldots, \si_{n-1},  a_1, \ldots, a_{2g}. $

  \end{itemize}

 \item Relations:

 \begin{itemize}

  \item[(R1)] $\si_i \si_j = \si_j \si_i$  \hspace{\stretch{1}}$(|i-j|\geq 2)$

  \item[(R2)] $\si_i \si_{i+1} \si_i = \si_{i+1} \si_i \si_{i+1}$
     \hspace{\stretch{1}}    $(1 \leq i \leq n-2)$ 

  \item[(R3)] $ a_1\cdots a_{2g}
        a_{1}^{-1}\cdots a_{2g}^{-1} =
        \si_{1}\cdots \si_{n-2} \si_{n-1}^2 \si_{n-2} \cdots \si_{1} $

  \item[(R4)] $  a_{r} A_{2,s}= A_{2,s} a_{r}$ \hspace{\stretch{1}}
    $(1\leq r,s \leq 2g; \;\;  r \neq s )$  

  \item[(R5)] $ \left( a_{1}\cdots a_{r}\right) A_{2,r}=\si_1^2 A_{2,r}
      \left( a_{1}\cdots a_{r}\right)$  \hspace{\stretch{1}} 
       $ (1\leq r\leq 2g) $
     
  \item[(R6)] $a_{r}\si_i=\si_i a_{r}$  \hspace{\stretch{1}} 
   $ (1\leq r \leq 2g; \; \; i\geq 2)$

  \end{itemize}

\end{itemize}
where 
$$
    A_{2,r}=\si_1^{-1}
   \left( a_{1}\cdots a_{r-1}a_{r+1}^{-1}\cdots a_{2g}^{-1}\right)
    \si_1^{-1}.
$$
\end{theorem}

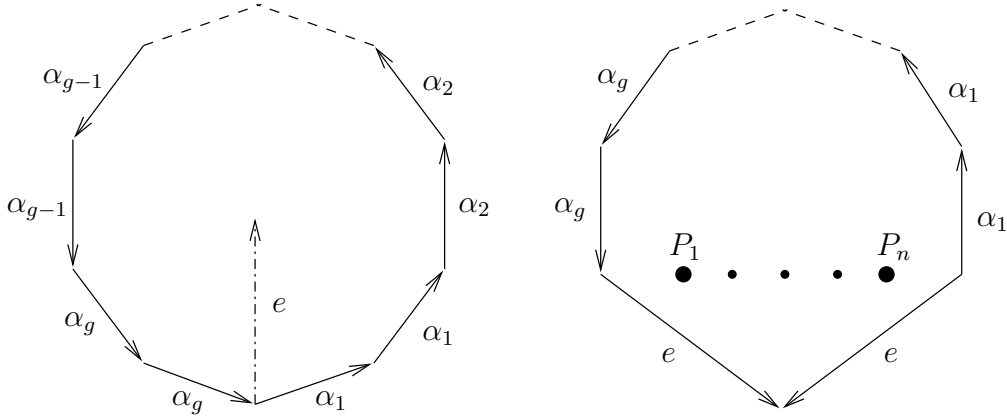
\begin{figure}[ht]
\centerline{\input{nonpolygon.pstex_t}\hspace{1cm}
              \input{nonpolygon2.pstex_t}}
\caption{Representation of a non-orientable surface $M$.}
\label{nonpolygon}
\end {figure}

\vspace{.6cm}
  Now we turn to the non-orientable case. Let $M$ be a closed non-orientable 
surface of genus $g\geq 2$. To 
represent a braid in $M$ we will also present 
the surface as a polygon, this time of $2g$ sides, as in Figure~\ref{nonpolygon},
and we make an additional cut: define the path $e$ as in the left hand side
of Figure~\ref{nonpolygon}, and cut the polygon along it. We  get $M$ 
represented as in the right hand side of the same figure, where we can 
also see how we choose the points $P_1,\ldots,P_n$. 

\begin{figure}[ht]
\centerline{\input{nona1r.pstex_t} \hspace{2cm} \input{nonsii.pstex_t}}
\caption{The generators of $B_n(M)$.}
\label{nona1r}
\end {figure}
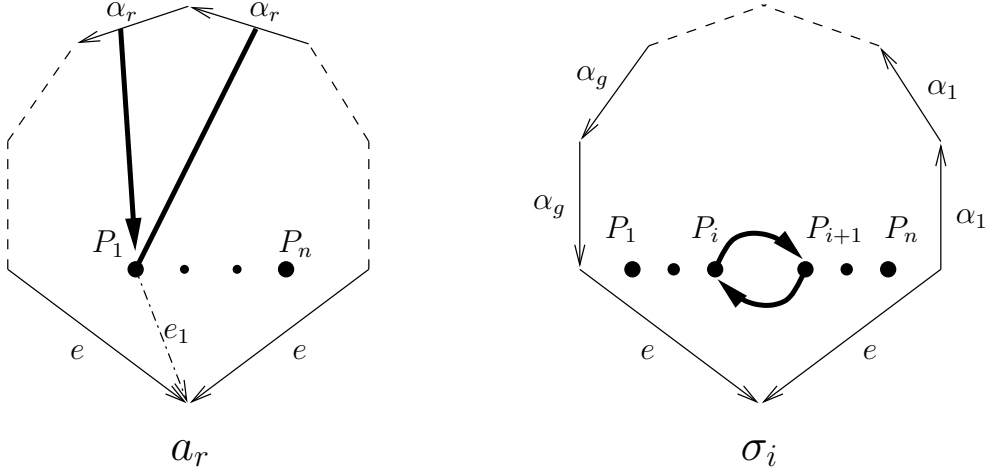

\vspace{.3cm}
  We define now the generators of $B_n(M)$. They will be similar 
to those of the orientable surface braid groups. For all 
$i\in\{1,\ldots n-1 \}$,
the braid $\si_i$ will be the same as in the orientable case.
For all $r\in\{ 1,\ldots,g\}$,
the braid $a_{r}$ consists on the first string passing through the
$r$-th wall, in the way of Figure~\ref{nona1r}, while the other strings are
trivial paths. 

\vspace{.3cm}
  There are six relations in the braid group of $M$ that are analogous to
those considered for an orientable surface. They can be shown to hold in 
the same way as in the orientable case; the only difference is the construction 
of the $P_1$-polygon. We denote by  $s_1,\ldots,s_g$ the first string of 
$a_1,\ldots,a_g$, respectively, and in this case we define another path, $e_1$, 
which goes from $P_1$ to the final point of $e$ (see Figure~\ref{nona1r}). Then
we cut along the paths $s_1,\ldots,s_g, e_1$ and glue along 
$\al_1,\ldots,\al_g, e$. The result is the {\it $P_1$-polygon} of $M$ 
whose sides, 
reading clockwise, are labeled by $s_1,s_1,s_2,s_2,\ldots, 
s_g, s_g, e_1, e_1^{-1}$.

\vspace{.3cm} 
 We claim that the six mentioned relations form a set of
defining relations of $B_n(M)$. To be more precise, we claim the following.

\vspace{.6cm}
\begin{theorem}\label{presBnN}
  If $M$ is a closed, non-orientable surface of genus $g\geq 2$, then 
$B_n(M)$ admits the following presentation:

\begin{itemize}

 \item Generators:
  \begin{itemize}

  \item[]   $\si_1, \ldots, \si_{n-1}, a_{1}, \ldots, a_{g}. $

  \end{itemize}

 \item Relations:

 \begin{itemize}

  \item[(r1)] $\si_i \si_j = \si_j \si_i$  \hspace{\stretch{1}}$(|i-j|\geq 2)$

  \item[(r2)] $\si_i \si_{i+1} \si_i = \si_{i+1} \si_i \si_{i+1}$
     \hspace{\stretch{1}}    $(1 \leq i \leq n-2)$ 

  \item[(r3)] $ a_{1}^2\cdots a_{g}^2=
        \si_{1}\cdots \si_{n-2} \si_{n-1}^2 \si_{n-2} \cdots \si_{1} $

  \item[(r4)] $  a_{r} A_{2,s}= A_{2,s} a_{r}$ \hspace{\stretch{1}}
    $(1\leq r,s \leq g; \;\;  r \neq s )$  

  \item[(r5)] $ \left( a_{1}^2\cdots a_{r-1}^2 a_{r}\right) A_{2,r}=
      \si_1^2 A_{2,r} \left( a_{1}^2\cdots a_{r-1}^2 a_{r}\right)$  
     \hspace{\stretch{1}}   $ (1\leq r\leq g) $
     
  \item[(r6)] $a_{r}\si_j=\si_ja_{r}$  \hspace{\stretch{1}} 
   $ (1\leq r \leq g; \; \; j\geq 2)$

  \end{itemize}

\end{itemize}
where 
$$
    A_{2,r}=\si_1^{-1}
   \left( a_{1}^2\cdots a_{r-1}^2 a_{r}^{-1}
      a_{r-1}^{-2}\cdots a_{1}^{-2}\right)
    \si_1.
$$
\end{theorem}

\section{A method for finding presentations}\label{method}

Consider an exact sequence of groups
$$
    1 \longrightarrow A \stackrel{i}{\longrightarrow} B  
     \stackrel{p}{\longrightarrow}   C  \longrightarrow 1, 
$$
where we suppose $A\subset B$, and $i$ is the inclusion map.
Suppose that $A$ and $C$ have presentations
$$
   A=<G_A; \: R_A>, \quad  \quad   C=<G_C; \: R_C>.
$$
 For each $y\in G_C$, we choose an element $\tilde{y}\in B$ such that 
$p(\tilde{y})=y$, and for each relator $r=y_1\ldots y_m \in R_C$, we write
$\tilde{r}=\tilde{y}_1\ldots \tilde{y}_m\in B$. Then it is clear that for 
every $r\in R_C$, there exists a word $f_r$ over $G_ A$ such that 
$\tilde{r}=f_r$ in $B$.

\vspace{.3cm}
 On the other hand, for all $x\in G_A$ and $y\in G_C$, there exists a
word $g_{x,y}$ over $G_ A$ such that 
$\tilde{y}\: x\: {\tilde{y}}^{-1}=g_{x,y}$ in $B$.

\begin{lemma}\label{baselem}
Under the above conditions,
$B$ admits the following presentation:

\begin{itemize}
\item Generators:$\quad  \{ G_A \} \cup \{ \tilde{y}; \; y\in G_C \}$

\item Relations:
\begin{itemize}
 \item {\bf Type 1:} $r_A=1$, \quad \quad \quad \quad for all $r_A\in R_A$.

 \item {\bf Type 2:} $\tilde{r}=f_r$, \quad \quad \quad \quad for all $r\in R_C$ .

 \item {\bf Type 3:} $\tilde{y}\: x\: {\tilde{y}}^{-1}=g_{x,y}$, 
       \quad for all $x\in G_A$, and all $y\in G_C$.
\end{itemize}
\end{itemize}
\end{lemma}

\vspace{.6cm}
  The proof of this lemma is left to the reader. The plan of the proofs of 
Theorems~\ref{presBnM} and \ref{presBnN} is as follows:

\vspace{.6cm}
\noindent {\bf Step 1.} We will introduce an abstract group 
$\overline{PB_n(M)}$ given by
its presentation, and define a homomorphism
$$
  \overline{PB_n(M)} \; \stackrel{\varphi}{\longrightarrow} \; PB_n(M).
$$

\vspace{.6cm}
\noindent {\bf Step 2.}  We will prove by induction on $n$ that $\varphi$ 
is an isomorphism, applying Lemma~\ref{baselem} to the exact 
sequence~(\ref{es2}):
$$
1\longrightarrow \pi_1(M\backslash {\cal P}', P_1) 
      \stackrel{u}{\longrightarrow} PB_n(M,{\cal P}) 
    \stackrel{v}{\longrightarrow} PB_{n-1}(M, {\cal P}') \longrightarrow 1.
$$

\vspace{.6cm}
\noindent {\bf Step 3.} 
We denote by $\overline{B_n(M)}$ the abstract group given by
the presentation of Theorem~\ref{presBnM} if $M$ is oriented, and by the
presentation of Theorem~\ref{presBnN} if $M$ is non-oriented. It is 
shown in Section~\ref{statements} that there is a well defined 
homomorphism 
$$
  \overline{B_n(M)} \; \stackrel{\psi}{\longrightarrow} \; B_n(M).
$$
We will apply Lemma~\ref{baselem} to the exact sequence~(\ref{es1}):
$$
    1\longrightarrow PB_n(M) \stackrel{e}{\longrightarrow}
       B_n(M) \stackrel{f}{\longrightarrow} \Si_n \longrightarrow 1
$$
to show that $\psi$ is actually an isomorphism.

\section{The braid groups of an orientable surface}\label{secori}
 In this section we prove Theorem~\ref{presBnM} following the procedure 
given in Section~\ref{method}. So, throughout the section, $M$ is 
assumed to be an orientable surface of genus $g\geq 1$.

\vspace{.6cm}
\noindent {\bf Step 1.} Let $\overline{PB_n(M)}$ be the group given by 
the following presentation: 
  
\vspace{.3cm}
\underline{\bf Presentation 1} 
\begin{itemize}

 \item Generators:  $\left\{a_{i,r}; \; \; 1\leq i \leq n, \; 1\leq r \leq 2g\right\}
    \cup  \left\{ T_{j,k}; \;\; 1\leq j < k \leq n \right\}$.

 \item Relations:

 \begin{itemize}

  \item[(PR1)] $ a_{n,1}^{-1}a_{n,2}^{-1}\cdots a_{n,2g}^{-1}
        a_{n,1}a_{n,2}\cdots a_{n,2g} =
        \prod_{i=1}^{n-1}{T_{i,n-1}^{-1}T_{i,n}}. $

  \item[(PR2)] $  a_{i,r} A_{j,s}= A_{j,s} a_{i,r} $ \hspace{\stretch{1}}
   $ (1\leq i < j \leq n; \: 1\leq r,s \leq 2g; \: r\neq s ).$  

   \item[(PR3)] $\left( a_{i,1}\cdots a_{i,r}\right) A_{j,r}
     \left( a_{i,r}^{-1}\cdots a_{i,1}^{-1} \right) 
     A_{j,r}^{-1}= T_{i,j} T_{i,j-1}^{-1}$ 
    \hspace{\stretch{1}} $(1\leq i < j \leq n; \:  1\leq r\leq 2g). $

  \item[(PR4)] $ T_{i,j} T_{k,l}=T_{k,l} T_{i,j} $   \hspace{\stretch{1}}
     $  (1 \leq i<j<k<l\leq n \; \mbox{ or } \;1\leq i<k<l \leq j \leq  n ).$

  \item[(PR5)] $T_{k,l} T_{i,j} T_{k,l}^{-1}=T_{i,k-1}T_{i,k}^{-1}T_{i,j}
         T_{i,l}^{-1} T_{i,k}T_{i,k-1}^{-1}T_{i,l}$ \hspace{\stretch{1}}
      $  (1\leq i<k \leq j<l\leq n).  $

  \item[(PR6)] $a_{i,r} T_{j,k}=T_{j,k} a_{i,r}$ \hspace{\stretch{1}} 
   $ (1 \leq i<j<k\leq n \; \mbox{ or } \; 1\leq j<k<i\leq n), \: 
    (1\leq r \leq 2g).$

  \item[(PR7)] $ a_{i,r} \left( a_{j,2g}^{-1}\cdots a_{j,1}^{-1} T_{j,k}
         a_{j,2g}\cdots a_{j,1}\right) = 
      \left( a_{j,2g}^{-1}\cdots a_{j,1}^{-1} T_{j,k}
         a_{j,2g}\cdots a_{j,1}\right) a_{i,r}$ 

\hspace{\stretch{1}}  $ (1\leq j<i\leq k\leq n).$

  \item[(PR8)] $ T_{j,n}=\left( \prod_{i=1}^{j-1}{a_{i,2g}^{-1}\cdots 
        a_{i,1}^{-1} T_{i,j-1} 
        T_{i,j}^{-1} a_{i,1}\cdots a_{i,2g}}\right) a_{j,1}\cdots 
        a_{j,2g}a_{j,1}^{-1}\cdots a_{j,2g}^{-1}.$

 \end{itemize}

\end{itemize}
Where 
$$
   A_{j,s}=a_{j,1}\cdots a_{j,s-1}a_{j,s+1}^{-1}\cdots a_{j,2g}^{-1}.
$$

\vspace{.3cm}
  Later, we will make use of a different presentation of $\overline{PB_n(M)}$,
based on the following lemma. 

\begin{lemma}\label{xlem}
 Let $F$ be the free group freely generated by $\{x_1,\ldots,x_{2g} \}$.
Set
$$
   X_r=x_1\cdots x_{r-1}x_{r+1}^{-1}\cdots x_{2g}^{-1}.  
$$
Then $\{X_1,\ldots X_{2g} \}$ is a free system of generators of $F$.
\end{lemma}

\begin{proof}
We only need to give the formulae of the change of generators, which are
$$
  x_k=\left( X_1 X_2^{-1}\cdots X_{k-2}X_{k-1}^{-1}\right)
     \left( X_{k+1}X_{k+2}^{-1}\cdots X_{2g-1}^{-1}X_{2g}\right)  
       \quad \quad
        \mbox{if $k$ is odd,}
$$
$$
  x_k^{-1}=\left( X_1 X_2^{-1}\cdots X_{k-2}^{-1}X_{k-1}\right) 
     \left( X_{k+1}^{-1}X_{k+2}\cdots X_{2g-1}^{-1}X_{2g} \right) 
        \quad \quad
        \mbox{if $k$ is even.}
$$
\end{proof}

\vspace{.3cm}
  As a direct consequence of this lemma,  $\overline{PB_n(M)}$ admits
the following presentation.

\vspace{.3cm}
\underline{\bf Presentation 2} 
\begin{itemize}

 \item Generators: $\left\{A_{i,r}; \; \; 1\leq i \leq n, \; 1\leq r \leq 2g\right\}
    \cup  \left\{ T_{j,k}; \;\; 1\leq j < k \leq n \right\}$.

 \item Relations: The same of Presentation 1, where
$$
  a_{i,k}=\left( A_{i,1} A_{i,2}^{-1}\cdots A_{i,k-2}A_{i,k-1}^{-1}\right)
     \left( A_{i,k+1}A_{i,k+2}^{-1}\cdots A_{i,2g-1}^{-1}A_{i,2g}\right)  
       \quad \quad
        \mbox{if $k$ is odd,}
$$
$$
  a_{i,k}^{-1}=\left( A_{i,1} A_{i,2}^{-1}\cdots A_{i,k-2}^{-1}A_{i,k-1}\right) 
     \left( A_{i,k+1}^{-1}A_{i,k+2}\cdots A_{i,2g-1}^{-1}A_{i,2g} \right) 
        \quad \quad
        \mbox{if $k$ is even.}
$$
\end{itemize}

\vspace{.3cm}
 According to  Step 1, we must define a homomorphism
$$
  \overline{PB_n(M)} \; \stackrel{\varphi}{\longrightarrow} \; PB_n(M).
$$
By abuse of notation, we will still denote by $a_{i,r}$ and $T_{i,j}$ the
braids that will be the images of $a_{i,r}$ and $T_{i,j}$, respectively, 
under the homomorphism $\varphi$. These braids are defined as follows.

\begin{figure}[ht]
\centerline{\input{air.pstex_t}\quad \input{ais.pstex_t}\quad 
  \input{Tij.pstex_t}}
\caption{The generators of $PB_n{M}$.}
\label{generators2}
\end {figure}
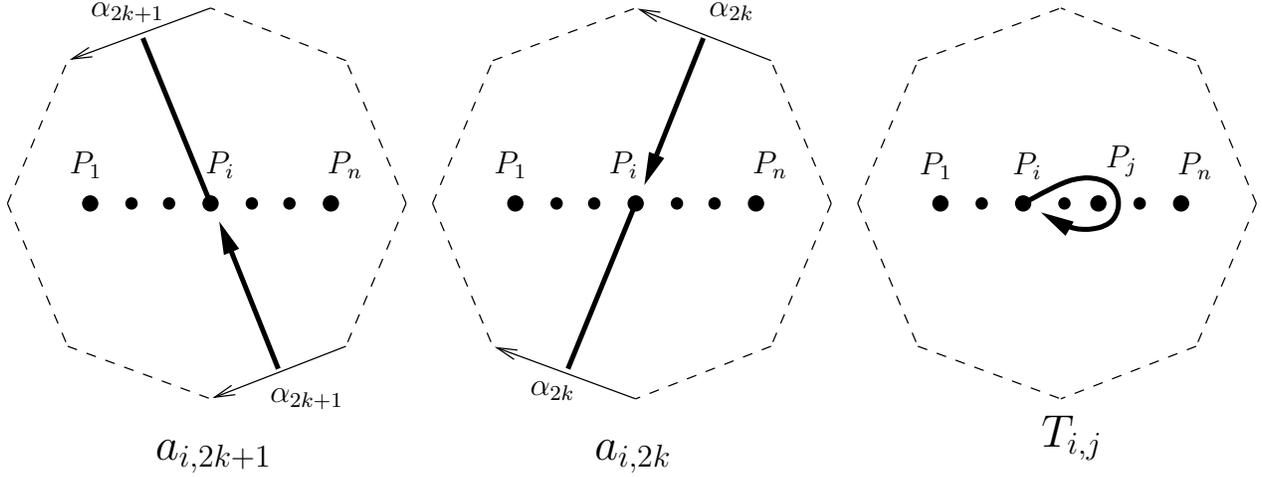

\begin{itemize}  
 \item In $a_{i,r}$, the $i$-th string goes through the $r$-th wall, as in
Figure~\ref{generators2}. This string will go upwards if $r$ is odd, and 
downwards otherwise. The other strings are trivial.
Note that $a_{1,r}=a_r$ for all $r$.

\item In $T_{i,j}$, the $i$-th string surrounds the points 
$P_{i+1},\ldots, P_{j}$, in the way of Figure~\ref{generators2}, while
the other strings are trivial paths. If $i=j$, 
we make $T_{i,j}$ to be the trivial braid.

\end{itemize}

\vspace{.3cm}
  We will denote by $s_{i,r}$ the $i$-th string of $a_{i,r}$,
and by $t_{i,j}$ that of $T_{i,j}$.
One can easily show that for any $i$, the set of paths 
$\{s_{i,1},\ldots,s_{i,2g}\}$ generates $\pi_1(M)$. Now, for any 
$i\in\{2,\ldots,n \}$ we can define the $P_i$-polygon as we defined the 
$P_1$-polygon in Section~\ref{statements}: we cut $L$ along 
$s_{i,1},\ldots,s_{i,2g}$ and glue along $\al_1,\ldots,\al_{2g}$.

\begin{figure}[ht]
\centerline{\input{Ajr.pstex_t}}
\caption{The braid $A_{j,r}$ in the $P_i$-polygon $(i<j)$.}
\label{Ajr}
\end {figure}
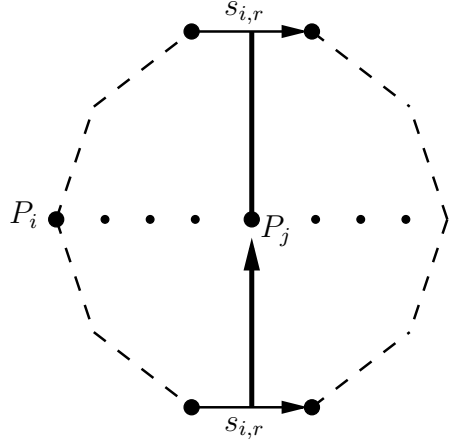

\vspace{.3cm}
 We define, for $2\leq j\leq n$ and $1\leq r \leq 2g$, the braid
$$
  A_{j,r}=a_{j,1}\cdots a_{j,r-1}a_{j,r+1}^{-1}\cdots a_{j,2g}^{-1}.
$$
Like in the representation of $A_{2,r}$ in the $P_1$-polygon considered
in Section~\ref{statements},  $A_{j,r}$ can be represented 
in the $P_i$-polygon (for $1\leq i < j$), as the braid of 
Figure~\ref{Ajr}, whose only 
nontrivial string is the $j$-th one, which goes upwards and crosses once 
the $r$-th wall $s_{i,r}$. Note that this representation does not depend on $i$,
but it is only valid when $i<j$.

\begin{figure}[ht]
\centerline{\input{rel81.pstex_t}\hspace{1.5cm} \input{rel82.pstex_t}}
\caption{The braid $a_{i,2g}^{-1}\cdots a_{i,1}^{-1}   T_{i,j-1}  
T_{i,j}^{-1} a_{i,1}\cdots a_{i,2g}$.}
\label{rel8}
\end {figure}
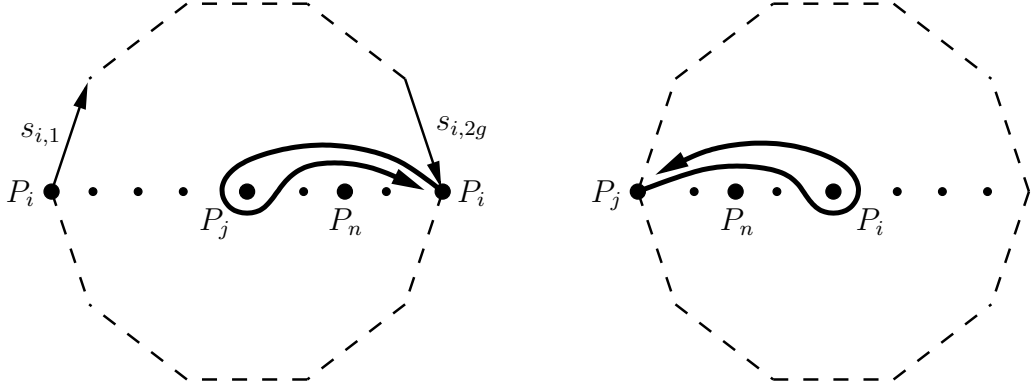

\vspace{.3cm}
 Now we define $\varphi$ in the obvious way. In order to show that it is a
homomorphism, we must show that the relations of $\overline{PB_n(M)}$  
still hold in $PB_n(M)$. Relations (PR4) and (PR5) can be easily 
checked, since they can be seen in the cylinder as if they were braids
over a disc (the interior of $L$). Relation (PR6) is obvious, once we have 
drawn the corresponding braids. Relations (PR1), (PR2) and (PR3) 
are analogous to Relations (R3), (R4) and (R5) of 
Theorem~\ref{presBnM}, and can be verified
in the same way. Relation (PR7) is easily checked in the $P_j$-polygon, and
finally, to verify Relation (PR8) we need all the $P_i$-polygons 
for $i=1,\ldots,j$: If $i<j$, it is clear by looking at the $P_i$-polygon
that
$$
a_{i,2g}^{-1}\cdots a_{i,1}^{-1}   T_{i,j-1}  T_{i,j}^{-1} 
a_{i,1}\cdots a_{i,2g}
$$
is equivalent to the braid on the left hand side of Figure~\ref{rel8}, 
thus it is equivalent to that on the right hand side, represented in the 
$P_j$-polygon. Then Relation (PR8) is clear, drawing all the factors
in the $P_j$-polygon.

\vspace{.3cm} 
 Hence, we have shown that $\varphi$ is a homomorphism, so this finishes 
the first step.

\vspace{.3cm}
 \noindent {\bf Step 2.} We show by induction on $n$ that 
$\varphi$ is an isomorphism. The case $n=1$  is clear, since the presentation of 
$\overline{PB_1(M)}$ turns to be
$$
 \overline{PB_1(M)}=\left< \{a_{1,1},\ldots,a_{1,2g}\}\; ; \; 
     a_{1,1}^{-1}a_{1,2}^{-1}\cdots a_{1,2g}^{-1} a_{1,1}a_{1,2}
     \cdots a_{1,2g} = 1  \right>,
$$
and this is also a presentation of $\pi_1(M)=PB_1(M)$. Moreover,
since $n=1$, one has $\varphi(a_{1,i})=a_{1,i}=s_{1,i}$ for all 
$i=1,\ldots,2g$, so $\overline{PB_1(M)}\stackrel{\varphi}{\simeq}PB_1(M)$.

\vspace{.3cm}
  Now suppose $\overline{PB_{n-1}(M)}\stackrel{\varphi}{\simeq}PB_{n-1}(M)$,
and recall the exact sequence~(\ref{es2}):
$$
    1\longrightarrow \pi_1(M\backslash {\cal P}', P_1) 
      \stackrel{u}{\longrightarrow} PB_n(M, {\cal P}) 
    \stackrel{v}{\longrightarrow} PB_{n-1}(M,{\cal P}') \longrightarrow 1.
$$

\vspace{.3cm}
 In order to apply Lemma~\ref{baselem} we need to know presentations of the 
groups at both hand sides. For the group on the left hand side, we have the 
presentation
$$
    \pi_1(M\backslash {\cal P}', P_1) =
    \left< \{s_{1,1},\ldots,s_{1,2g},t_{1,2},\ldots,t_{1,n-1}\}; \; \; 
     \phi \; \; \right>.
$$

\vspace{.3cm}
   It will be good for our purposes to include $t_{1,n}$ among the generators, 
so we add a single relation which can be easily deduced from the pictures
(using the $P_1$-polygon):
$$
 \pi_1(M\backslash {\cal P}', P_1) =
    \left< \{s_{1,1},\ldots,s_{1,2g},t_{1,2},\ldots,t_{1,n}\}; \; \; 
     t_{1,n}=  s_{1,1}\cdots s_{1,2g}
     s_{1,1}^{-1}\cdots s_{1,2g}^{-1} \; \; \right>.
$$

\vspace{.3cm}
  We know as well, by the induction hypothesis, two presentations of 
$PB_{n-1}(M)$; we shall use Presentation 2 of $\overline{PB_{n-1}(M)}$. 
So we can apply Lemma~\ref{baselem} to the exact sequence~(\ref{es2}).

\vspace{.3cm}
Note that $v(a_{i,r})= a_{i-1,r}$, for 
$i=2,\ldots, n$, 
so $v(A_{i,r})= A_{i-1,r}$, for $i=2,\ldots, n$. 
Note also that $v(T_{i,j})= T_{i-1,j-1}$, where 
$2\leq i\leq j \leq n$. So we know pre-images by $v$
of the generators of $PB_{n-1}(M,{\cal P}')$.  
 
\vspace{.3cm}
  It is also clear that $u(s_{1,r})=a_{1,r}$ and $u(t_{1,j})=T_{1,j}$ for 
all possible $r$ and $j$. Hence, we obtain 
immediately that a set of generators of $PB_n(M,{\cal P})$ is 
$$
   \{ a_{1,r}; \;\; 1\leq r \leq 2g \} \cup 
   \{ A_{i,r}; \;\; 2\leq i \leq n, \; 1\leq r \leq 2g\} \cup
   \{ T_{j,k}; \;\; 1\leq j < k \leq n \}.
$$
We can apply again Lemma~\ref{xlem} to have a new set of generators
$$
\{ a_{i,r}; \;\; 1\leq i \leq n, \; 1\leq r \leq 2g \} \cup 
   \{ T_{j,k}; \;\; 1\leq j < k \leq n \}.
$$
which is the image by $\varphi$  of the generating set of 
$\overline{PB_n(M)}$. In particular, $\varphi$ is surjective. 

\vspace{.3cm}
  Now we prove that $\varphi$ is an isomorphism by the following procedure.

\vspace{.3cm}
  First, we denote  by $G_A$ the set of generators of 
$\pi_1(M\backslash {\cal P}', P_1)$, and by $G$ the set of generators 
of $\overline{PB_n(M)}$. We consider the unique relation in the presentation 
of $\pi_1(M\backslash {\cal P}', P_1)$, which we can consider via $u$ as a 
relation in $PB_n(M)$. This will be the unique relation of Type 1 in the 
presentation of $PB_n(M)$. The procedure starts by showing that this relation
holds when it is considered in $\overline{PB_n(M)}$, that is, we have a relation
in $\overline{PB_n(M)}$ which maps by $\varphi$ to the only relation in the 
presentation of $\pi_1(M\backslash {\cal P}', P_1)$.

\vspace{.3cm}
  Next, for each relator $r$ of $PB_{n-1}(M)$, we
consider the ``canonical'' pre-image by $v$ of $r$, denoted by $\tilde{r}$, 
in the way of Lemma~\ref{baselem}. Since $PB_n(M)$ and $\overline{PB_n(M)}$
have the ``same'' generators (via $\varphi$), we can also consider 
$\tilde{r}$ as a word over $G$. Now we find a word $U$ over $G$ such that the 
equality $\tilde{r}=U$ holds in $\overline{PB_n(M)}$, and such that $\varphi(U)$ 
is a word over $G_A$. This will give us the relations of Type 2 in the 
presentation of $PB_n(M)$.

\vspace{.3cm}
  Finally, for each $x\in G_A$ and each generator $y$ of $PB_{n-1}(M)$, we 
find a word $V$ over $G$ such that the equality 
$\tilde{y}\: x \: \tilde{y}^{-1}=V$ holds in $\overline{PB_n(M)}$, 
where $\tilde{y}$ is the canonical pre-image by $v$ of $y$, and such that 
$\varphi(V)$ is a word over $G_A$. This will give us the relations of Type 3 
in the presentation of $PB_n(M)$.

\vspace{.3cm}
In this way, we will have found all relations
of Types 1, 2 and 3 of Lemma~\ref{baselem} and, therefore, a presentation
of $PB_n(M)$, and, at the same time, we 
will have shown that $\varphi$ is injective, and 
consequently, that $\varphi$ is an isomorphism. 

\vspace{.3cm} 
  Let us start with the procedure. The unique relation in the presentation
of $\pi_1(M\backslash {\cal P}', P_1)$ corresponds to Relation (PR8) of  
$\overline{PB_n(M)}$, for $j=1$, so it holds in this group. 

\vspace{.3cm}
  Relations of Type $2$ are easy to find. First, (PR1) can be seen as follows:
$$
    a_{n,1}^{-1}a_{n,2}^{-1}\cdots a_{n,2g}^{-1}
        a_{n,1}a_{n,2}\cdots a_{n,2g} \left(
        \prod_{i=2}^{n-1}{T_{i,n-1}^{-1}T_{i,n}}\right)^{-1}
   = T_{1,n-1}^{-1} T_{1,n}. 
$$
Note that the left hand side maps by $\varphi$ to $\tilde{r}$, where
$r$ is a relator of $PB_{n-1}(M)$ corresponding to (PR1), while the 
right hand side maps by $\varphi$ to a word over $G_A$. Hence, $U$ is
equal to the right hand side of the equation, 
and this yields the first relation of Type~2. The remaining relations 
of Type 2 are also images by $\varphi$ of relations in Presentation~1; 
namely (PR2), (PR3), (PR4) and (PR5) when $i\geq 2$, (PR6) when $i\geq 2$
and $j\geq 2$, and (PR7), (PR8) when $j\geq 2$. For all these relations, the 
word $U$ is just the trivial word, except for (PR8), for which 
$U=  a_{1,2g}^{-1} \cdots a_{1,1}^{-1} T_{1,j-1} T_{1,j}^{-1} a_{1,1}
\cdots a_{1,2g}$.

\vspace{.3cm} 
  Finally, we find the relations of  Type $3$. For $i=1$, (PR2) becomes
$$
    A_{j,s} a_{1,r} A_{j,s}^{-1} = a_{1,r}  \quad \quad \quad (r\neq s),
$$
so $V=a_{1,r}$. Next, using (PR2), Relation (PR3) turns to be equivalent to
$$
  A_{j,r} a_{1,r} A_{j,r}^{-1} =
  \left( a_{1,r-1}^{-1}\cdots a_{1,1}^{-1} \right)
   T_{1,j-1} T_{1,j}^{-1}
  \left( a_{1,1}\cdots a_{1,r}\right),
$$
so $V$ equals the right hand side of this equation. Relations of the form 
$T_{k,l} T_{1,j}T_{k,l}^{-1} = V$, where $V$ is a word over $G_A$, 
follow from (PR4)-(PR5),
while those of the form  $T_{k,l} a_{1,r} T_{k,l}^{-1} = V$
follow from (PR6), when $i=1$. Also, if $j>k$, we obtain from (PR6) the relations 
$A_{j,r}T_{1,k}A_{j,r}^{-1}= V$, where $V$ is a word over $G_A$.

\vspace{.3cm}
 The only remaining relations are those of the form $A_{j,r}T_{1,k}A_{j,r}^{-1}=
V$, when  $1 < j \leq k$, which are deduced as follows: By (PR7), we know that 
$a_{j,s}$ commutes with the element 
$$
  a_{1,2g}^{-1}\cdots a_{1,1}^{-1} T_{1,k}
         a_{1,2g}\cdots a_{1,1}
$$
for $s=1,\ldots,2g$. This implies that $A_{j,r}$ commutes with the
same element, so
$$
  \left(a_{1,2g}^{-1}\cdots a_{1,1}^{-1} T_{1,k}
         a_{1,2g}\cdots a_{1,1}\right) =
  A_{j,r}\left(a_{1,2g}^{-1}\cdots a_{1,1}^{-1} T_{1,k}
         a_{1,2g}\cdots a_{1,1}\right) A_{j,r}^{-1}
$$
$$
 = \left( A_{j,r}a_{1,2g}^{-1}A_{j,r}^{-1}\right) \cdots 
   \left( A_{j,r} a_{1,1}^{-1} A_{j,r}^{-1}\right) 
   \left( A_{j,r} T_{1,k} A_{j,r}^{-1} \right) 
   \left( A_{j,r} a_{1,2g} A_{j,r}^{-1}\right) \cdots 
   \left( A_{j,r} a_{1,1} A_{j,r}^{-1} \right). 
$$
But using (PR2) and (PR3) 
we know how to write all the terms in the above product (except the
middle one) as words over $G_A$, so we are done.

\vspace{.3cm}
  Hence, we have shown that 
$\overline{PB_n(M)}\stackrel{\varphi}{\simeq}PB_n(M)$ and therefore, we 
have proved:

\vspace{.3cm}
\begin{theorem}\label{presPBnM}
  If $M$ is a closed, orientable surface of genus $g\geq 1$, then 
$PB_n(M)$ admits Presentation 1 (and also Presentation 2) as presentation.
\end{theorem}

\vspace{.6cm}
\noindent {\bf Step 3.} Now we want to find a presentation of $B_n(M)$, for 
$g\geq 1$. We define then the group $\overline{B_n(M)}$, given by the 
presentation in Theorem~\ref{presBnM}. 

\vspace{.3cm}
  This is the most reduced presentation we have found. But to show its
validity we will modify it, obtaining a new one with 
more generators and relations, but equivalent to the first one.

\vspace{.3cm}
  First, we change our notation, and call $a_{1,r}$ the generators $a_r$, for 
$r=1,\ldots,2g$. Then we must simply add to the given presentation the generators

\begin{itemize}
 \item[$\relbar$] $ a_{i,r}  \quad \quad \quad i=2,\ldots,n\;; 
                           \quad r=1,\ldots,2g, $

 \item[$\relbar$] $ T_{j,k}  \quad \quad \quad  1\leq j < k \leq n,  $
\end{itemize}
and the relations
\begin{itemize}
 \item[(R7)] $a_{j+1,r}=\si_j a_{j,r}\si_j $   \hspace{\stretch{1}}
     ( $1\leq j \leq n-1; \; 1\leq r \leq 2g; \; r$ even).

 \item[(R8)] $a_{j+1,r}=\si_j^{-1} a_{j,r}\si_j^{-1} $   \hspace{\stretch{1}}
      ($1\leq j \leq n-1; \; 1\leq r \leq 2g; \; r$ odd).

 \item[(R9)]  $ T_{j,k}= \si_{j} \si_{j+1}\cdots \si_{k-2} 
          \si_{k-1}^2  \si_{k-2} \cdots \si_j $ \hspace{\stretch{1}}
      $(1 \leq j < k \leq n)$.

\end{itemize}

\vspace{.3cm}
  Clearly, both presentations define the same group, that is,
$\overline{B_n(M)}$.  Now we define 
$\psi:\:\overline{B_n(M)}\rightarrow B_n(M)$ in the natural way.
It is an easy exercise to show, using the same methods as before, that 
Relations (R7), (R8) and (R9) map to relations in $B_n(M)$. 
Therefore, $\psi$ is a well defined homomorphism.

\vspace{.3cm}
   Recall now the exact sequence~(\ref{es1}):
$$
   1\longrightarrow PB_n(M) \stackrel{e}{\longrightarrow}
       B_n(M) \stackrel{f}{\longrightarrow} \Si_n \longrightarrow 1.
$$ 

\vspace{.3cm}
  We know by Theorem~\ref{presPBnM} a presentation of $PB_n(M)$ (say
Presentation 1), and 
it is also known that a presentation of $\Si_n$ is
\begin{itemize} 
  \item Generators: $ \de_1,\ldots, \de_{n-1}.$     

  \item Relations:
     \begin{itemize}

        \item $\de_i \de_j = \de_j \de_i$ \quad \quad \quad \quad 
                $|i-j|\geq 2$,

	\item $\de_i \de_{i+1} \de_i = \de_{i+1} \de_i \de_{i+1}$   
              \quad  \quad 
              $1\leq i \leq n-2,$ 

	\item $\de_i^2=1$ \quad \quad \quad \quad \quad \quad \quad 
              $1\leq i \leq n-1,$
     \end{itemize}
\end{itemize}
where $\de_i$ is the permutation $(i, \; i+1 )$, for any $i$. 

\vspace{.3cm}
  Now $\si_i$ is clearly a  pre-image by $f$ of $\de_i$, so by 
Lemma~\ref{baselem} $\overline{B_n(M)}$ and $B_n(M)$ have the same 
generators, and $\psi$ is surjective.

\vspace{.3cm}
  Similarly to what we did in Step 2, we show now that $\psi$ is an 
isomorphism by the following procedure.

\vspace{.3cm}
  First, we denote by $G_A$ the set of generators of $PB_n(M)$, and by $G$ 
the set of generators of $\overline{B_n(M)}$. For each relation in the 
presentation of $PB_n(M)$, we consider it via $e$ as a relation in $B_n(M)$, 
and we show that it also holds in $\overline{B_n(M)}$. 

\vspace{.3cm}
  Next, for each relator $r$ of $\Si_n$, we consider its canonical pre-image by
$f$, denoted by $\tilde{r}$. Then we find a word $U$ over $G$
such that the equality $\tilde{r}=U$ holds in $\overline{B_n(M)}$, and such that 
$\psi(U)$ is a word over $G_A$.

\vspace{.3cm}
  Finally, for each $x\in G_A$ and each generator $\de_i$ of $\Si_n$, we find a 
word $V$ over $G$ such that the equality $\si_i \: x \: \si_i^{-1} = V$ holds
in $\overline{B_n(M)}$, and such that $\psi(V)$ is a word over $G_A$.

\vspace{.3cm}
  This gives us the relations of Types 1, 2, and 3 of Lemma~\ref{baselem} 
and, therefore, a presentation of $B_n(M)$, and,
at the same time, this shows that $\psi$ is injective, and, consequently,
that $\psi$ is an isomorphism.

\vspace{.3cm}
   Let us then verify in $\overline{B_n(M)}$ the relations of Type $1$. In the 
case of (PR1), we start with (R3):
\begin{equation}\label{vuelta}
    a_{1,1}\cdots a_{1,2g} a_{1,1}^{-1}\cdots a_{1,2g}^{-1} =
        \si_{1}\cdots \si_{n-2} \si_{n-1}^2 \si_{n-2} \cdots \si_{1}. 
\end{equation}
Using (R7) and (R8), we see that
the left hand side of Equation~(\ref{vuelta}) becomes
$$
  \si_{1}\cdots \si_{n-1} \left( a_{n,1}\cdots a_{n,2g}\right)   
   \si_{n-1}^{-1} \cdots \si_{1}^{-2} \cdots \si_{n-1}^{-1} 
  \left( a_{n,1}^{-1}\cdots a_{n,2g}^{-1} \right) \si_{n-1} \cdots \si_{1}.
$$
On the other hand, from (R1), (R2) (braid relations) and 
(R9), we get
$$
   T_{i,n-1}^{-1}T_{i,n}\; = \; \si_i^{-1} \si_{i+1}^{-1} \cdots \si_{n-2}^{-1}
    \si_{n-1}^2 \si_{n-2} \cdots \si_i \; = \;
    \si_{n-1} \cdots \si_{i+1} \si_i^{2} \si_{i+1}^{-1} \cdots
    \si_{n-1}^{-1},
$$
so
$$
   \prod_{i=1}^{n-1}{T_{i,n-1}^{-1}T_{i,n}}=\si_{n-1} \cdots \si_1^{2} 
    \cdots \si_{n-1}.
$$
Therefore, Equation~(\ref{vuelta}) becomes 
$$
   a_{n,1}\cdots a_{n,2g}   
  \left( \prod_{i=1}^{n-1}{T_{i,n-1}^{-1}T_{i,n}}\right)^{-1}
   a_{n,1}^{-1}\cdots a_{n,2g}^{-1}=1,
$$   
which is clearly equivalent to (PR1).

\vspace{.3cm}
We will use in what follows some relations of $\overline{B_n(M)}$
easily deduced from (R1)-(R9). From (R7) and (R8), we get
\begin{equation}
  a_{i,r}= \left( \si_{i-1}^{-1}\cdots \si_{1}^{-1} \right) a_{1,r}
           \left( \si_{1}^{-1}\cdots \si_{i-1}^{-1} \right)
    \quad \mbox{if $r$ is odd.}
\end{equation}
\begin{equation}
  a_{i,r}= \left( \si_{i-1}\cdots \si_{1} \right) a_{1,r}
           \left( \si_{1}\cdots \si_{i-1} \right)
    \quad \mbox{if $r$ is even.}
\end{equation}
\begin{equation}
  A_{j,s}= \left( \si_{j-1}^{-1}\cdots \si_{2}^{-1} \right) A_{2,s}
           \left( \si_{2}^{-1}\cdots \si_{j-1}^{-1} \right)=
           \left( \si_{j-1}^{-1}\cdots \si_{1}^{-1} \right) A_{1,s}
           \left( \si_{1}^{-1}\cdots \si_{j-1}^{-1} \right).
\end{equation}

\vspace{.3cm}
  Also, from (R1) and (R2), we obtain
\begin{equation}
 \si_{j} \left( \si_k \si_{k-1}\cdots \si_{i} \right) =
         \left( \si_k \si_{k-1}\cdots \si_{i} \right) \si_{j+1}
  \quad \quad (i\leq j <k).
\end{equation}
\begin{equation}
 \si_{j} \left( \si_{k}^{-1}\si_{k-1}^{-1}\cdots \si_{i}^{-1} \right)=
	\left( \si_{k}^{-1}\si_{k-1}^{-1}\cdots \si_{i}^{-1} \right)
  \si_{j+1}  \quad \quad (i\leq j <k).
\end{equation}
\begin{equation}
 \si_{i}\cdots \si_{k-1} \si_{k}^{2}\si_{k-1}^{-1}\cdots \si_{i}^{-1} =
\si_{k}^{-1}\cdots \si_{i+1}^{-1}\si_{i}^{2}\si_{i+1}\cdots \si_{k}.
\end{equation}

\vspace{.3cm}
  Now using (6), (7), (8) and (R6), we see that if $1\leq k\leq j-2$;
$$
\begin{array}{l}
  \si_k A_{j,s} = \si_{k} \left( \si_{j-1}^{-1}\cdots \si_{1}^{-1} \right)
                 A_{1,s} \left(\si_{1}^{-1}\cdots \si_{j-1}^{-1}\right) \quad \quad \\
     \\
   =\left( \si_{j-1}^{-1}\cdots \si_{1}^{-1} \right) \si_{k+1}
          A_{1,s} \left(\si_{1}^{-1}\cdots \si_{j-1}^{-1}\right)  
\end{array}
$$
\begin{equation}
  =\left( \si_{j-1}^{-1}\cdots \si_{1}^{-1} \right) 
     A_{1,s} \si_{k+1}\left(\si_{1}^{-1}\cdots \si_{j-1}^{-1}\right) =
    A_{j,s} \si_k.
\end{equation}
In the same way, using (6), (R6) and (R4), we get
$$
 a_{1,r} A_{j,s}= a_{1,r}\left( \si_{j-1}^{-1}\cdots \si_{2}^{-1} \right)
   A_{2,s} \left(\si_{2}^{-1}\cdots \si_{j-1}^{-1}\right)=
    A_{j,s} a_{1,r},
$$
if $r\neq s$ and $1<j$. 

\vspace{.3cm}
  Therefore, if $i<j$ and $r\neq s$, by (4) and (5) $a_{i,r}$ is a product
of elements which commute with $A_{j,s}$, so we obtain
$$
   a_{i,r} A_{j,s}= A_{j,s} a_{i,r},
$$
which shows that (PR2) holds in $\overline{B_n(M)}$.

\vspace{.3cm}
  Now we verify Relation (PR3). We will do the case when
$r$ is odd, the other case being analogous. It is clear which of the
known relations of $\overline{B_n(M)}$ we are using at each of the 
following equalities:
$$
\begin{array}{l}
  \left( a_{i,1}\ldots a_{i,r} \right) A_{j,r} =
  \left( \si_{i-1}^{-1}\cdots \si_{1}^{-1} \right)
  \left( a_{1,1}\ldots a_{1,r} \right) 
  \left(\si_{1}^{-1}\cdots \si_{i-1}^{-1}\right) A_{j,r} \\ \\

 =  \left( \si_{i-1}^{-1}\cdots \si_{1}^{-1} \right)
  \left( a_{1,1}\ldots a_{1,r} \right) A_{j,r}
 \left(\si_{1}^{-1}\cdots \si_{i-1}^{-1}\right) \\ \\

 = \left( \si_{i-1}^{-1}\cdots \si_{1}^{-1} \right)
 \left( a_{1,1}\ldots a_{1,r} \right) 
  \left( \si_{j-1}^{-1}\cdots \si_{2}^{-1} \right)
   A_{2,r} \left(\si_{2}^{-1}\cdots \si_{j-1}^{-1}\right)
 \left(\si_{1}^{-1}\cdots \si_{i-1}^{-1}\right) \\ \\

=  \left( \si_{i-1}^{-1}\cdots \si_{1}^{-1} \right)
\left( \si_{j-1}^{-1}\cdots \si_{2}^{-1} \right) 
\left( a_{1,1}\ldots a_{1,r} \right)  A_{2,r}
\left(\si_{2}^{-1}\cdots \si_{j-1}^{-1}\right)
 \left(\si_{1}^{-1}\cdots \si_{i-1}^{-1}\right) \\ \\

= \left( \si_{i-1}^{-1}\cdots \si_{1}^{-1} \right)
\left( \si_{j-1}^{-1}\cdots \si_{2}^{-1} \right) \si_1^2
 A_{2,r}\left( a_{1,1}\ldots a_{1,r} \right)
\left(\si_{2}^{-1}\cdots \si_{j-1}^{-1}\right)
 \left(\si_{1}^{-1}\cdots \si_{i-1}^{-1}\right) \\ \\

 = \left(\si_{i}\cdots \si_{j-2} \si_{j-1}^{2}  \si_{j-2}^{-1}\cdots
     \si_{1}^{-1}\right)
\left( \si_{j-1}^{-1}\cdots \si_{2}^{-1}  A_{2,r} 
\si_{2}^{-1}\cdots \si_{j-1}^{-1}\right)
\left( a_{1,1}\ldots a_{1,r} \right) 
\left(\si_{1}^{-1}\cdots \si_{i-1}^{-1}\right) \\ \\

 = \left(\si_{i}\cdots \si_{j-2} \si_{j-1}^{2}  \si_{j-2}^{-1}\cdots
     \si_{1}^{-1}\right)    A_{j,r} 
\left( a_{1,1}\ldots a_{1,r} \right) 
\left(\si_{1}^{-1}\cdots \si_{i-1}^{-1}\right) \\ \\

 = \left(\si_{i}\cdots \si_{j-2} \si_{j-1}^{2}  \si_{j-2}^{-1}\cdots
     \si_{i}^{-1}\right)    A_{j,r} 
\left( a_{i,1}\ldots a_{i,r} \right) \\ \\

= T_{i,j} T_{i,j-1}^{-1} A_{j,r}  \left( a_{i,1}\ldots a_{i,r} \right).
\end{array}
$$

\vspace{.3cm}
  This shows the case of (PR3). Relations (PR4) and (PR5) are actually
relations in the braid group of the disc, so they are a consequence of
(R1) and (R2). (PR6) is obtained easily from (R9), (4), (5) and the
braid relations (R1) and (R2). So we may turn to (PR7): It is clear 
that it suffices to show that in $\overline{B_n(M)}$, $A_{i,r}$
commutes with $\left( a_{j,2g}^{-1}\cdots a_{j,1}^{-1} T_{j,k}
a_{j,2g}\cdots a_{j,1}\right)$ for $1\leq j < i \leq k < n$. 
This is shown as follows (remember that we can already use 
(PB1)-(PB6)):
$$
\begin{array}{l}
   A_{i,r} \left( a_{j,2g}^{-1}\cdots a_{j,1}^{-1} T_{j,k}
   a_{j,2g}\cdots a_{j,1}\right) \\ \\

 = \left( a_{j,2g}^{-1}\cdots a_{j,r+1}^{-1} \right) A_{i,r}
   \left( a_{j,r}^{-1}\cdots a_{j,1}^{-1} \right) T_{j,k}
   a_{j,2g}\cdots a_{j,1} \\ \\

 = \left( a_{j,2g}^{-1}\cdots a_{j,1}^{-1} \right) 
   T_{j,i} T_{j,i-1}^{-1} A_{i,r} T_{j,k}  a_{j,2g}\cdots a_{j,1} \\ \\

 = \left( a_{j,2g}^{-1}\cdots a_{j,1}^{-1} \right) 
   T_{j,i} T_{j,i-1}^{-1} A_{i,r} 
   \left(  \si_j \cdots \si_{k-1}^2 \cdots \si_j \right)  
   a_{j,2g}\cdots a_{j,1} \\ \\

 = \left( a_{j,2g}^{-1}\cdots a_{j,1}^{-1} \right) 
   T_{j,i} T_{j,i-1}^{-1} A_{i,r} 
   \left(  \si_j \cdots \si_{k-1}^2 \cdots \si_1 \right)  
   a_{1,2g}\cdots a_{1,1} 
   \left( \si_1^{-1} \cdots \si_{j-1}^{-1} \right) \\ \\

 \mbox{(using (R3))}  \\ \\

 = \left( a_{j,2g}^{-1}\cdots a_{j,1}^{-1} \right) 
   T_{j,i} T_{j,i-1}^{-1} A_{i,r} 
   \left(  \si_j \cdots \si_{k-1} \si_{k}^{-1} \cdots \si_{n-1}^{-2}
   \cdots \si_1^{-1} \right)  
   a_{1,1}\cdots a_{1,2g} 
   \left( \si_1^{-1} \cdots \si_{j-1}^{-1} \right) \\ \\

\mbox{(by (R9) and (10))} \\ \\

 =  a_{j,2g}^{-1}\cdots a_{j,1}^{-1} 
   \left( \si_j \cdots \si_{i-1}^2  A_{i,r} 
   \si_{i-1} \cdots \si_{k-1} \si_{k}^{-1} \cdots 
   \si_{n-1}^{-2} \cdots \si_1^{-1} \right)  
   a_{1,1}\cdots a_{1,2g} 
   \left( \si_1^{-1} \cdots \si_{j-1}^{-1} \right) \\ \\

 =  a_{j,2g}^{-1}\cdots a_{j,1}^{-1} 
   \left( \si_j \cdots \si_{i-1}  A_{i-1,r} 
   \si_{i} \cdots \si_{k-1} \si_{k}^{-1} \cdots 
   \si_{n-1}^{-2} \cdots \si_1^{-1} \right)  
   a_{1,1}\cdots a_{1,2g} 
   \left( \si_1^{-1} \cdots \si_{j-1}^{-1} \right) \\ \\

 =  a_{j,2g}^{-1}\cdots a_{j,1}^{-1} 
   \left( \si_j \cdots \si_{k-1} \si_k^{-1} \cdots \si_{n-1}^{-2}
   \cdots \si_i^{-1}\si_{i-1}\si_{i-2}^{-1} \cdots \si_1^{-1} \right) 
   A_{i,r} a_{1,1}\cdots a_{1,2g}
   \si_1^{-1} \cdots \si_{j-1}^{-1}  \\ \\

 =  a_{j,2g}^{-1}\cdots a_{j,1}^{-1} 
   \left( \si_j \cdots \si_{k-1} \si_k^{-1} \cdots \si_{n-1}^{-2}
   \cdots \si_1^{-1} \right)  a_{1,1}\cdots a_{1,2g}
   \left( \si_1^{-1} \cdots \si_{j-1}^{-1} \right) A_{i,r} \\ \\

\mbox{(by (R3) again)} \\ \\

 = a_{j,2g}^{-1}\cdots a_{j,1}^{-1} T_{j,k}
  \left( \si_{j-1}\cdots \si_1 a_{1,2g}\cdots a_{1,1} 
  \si_1^{-1} \cdots \si_{j-1}^{-1} \right) A_{i,r} \\ \\

 = \left( a_{j,2g}^{-1}\cdots a_{j,1}^{-1} T_{j,k}
  a_{j,2g}\cdots a_{j,1} \right) A_{i,r}.
\end{array}
$$

\vspace{.3cm}
 Finally, Relation (PR8) is verified using some intermediary results.
The first is evident: by (R4) we see that in $\overline{B_n(M)}$, 
$A_{1,2g} A_{2,2g} = A_{2,2g} A_{1,2g}$,
and moreover this braid commutes with $\si_1$, since
$$
   A_{1,2g} A_{2,2g} \si_1 =  A_{1,2g} \si_1^{-1} A_{1,2g} =
  \si_1 A_{2,2g} A_{1,2g} =  \si_1 A_{1,2g} A_{2,2g}.
$$
Analogously, one shows that  $(a_{1,2g} a_{2,2g})$ commutes with $\si_1$.
The following result is a consequence of the previous ones and of (R5): 
$$
\begin{array}{l}
   a_{1,2g} A_{2,2g} a_{1,2g}^{-1} = 
    \left( a_{1,2g-1}^{-1} \cdots  a_{1,1}^{-1}\right) \si_1^2 A_{2,2g}
    \left( a_{1,1} \cdots a_{1,2g-1} \right) \\ \\

  =  A_{1,2g}^{-1}  \si_1^2 A_{2,2g} A_{1,2g} = 
     A_{1,2g}^{-1}   A_{2,2g} A_{1,2g} \si_1^2 =
     A_{2,2g} \si_1^2,
\end{array}
$$
so we obtain
\begin{equation}
   a_{1,2g}^{-1} A_{2,2g} = A_{2,2g} a_{1,2g}^{-1} \si_1^{-2}.
\end{equation}

\vspace{.3cm} 
  Now we consider the factors in the right hand side of (PR8), and we
see that 
$$
\begin{array}{l}
 \left( a_{i,2g}^{-1}\cdots  a_{i,1}^{-1}\right)  T_{i,j-1} 
  T_{i,j}^{-1} \left( a_{i,1}\cdots a_{i,2g} \right) \\ \\

 =  \left( a_{i,2g}^{-1}\cdots  a_{i,1}^{-1}\right)
  \si_i \cdots \si_{j-2} \si_{j-1}^2 \si_{j-2}^{-1} \cdots \si_{i}^{-1}
    \left( a_{i,1}\cdots a_{i,2g} \right) \\ \\

 = \si_{i-1}^{-1} \cdots \si_{1}^{-1} 
   \left( a_{1,2g}^{-1}\cdots  a_{1,1}^{-1}\right)
   \si_1 \cdots \si_{j-2} \si_{j-1}^2 \si_{j-2}^{-1} \cdots \si_{1}^{-1} 
   \left( a_{1,1}\cdots a_{1,2g} \right)
   \si_1 \cdots \si_{i-1} \\ \\

\mbox{(by (9))} \\ \\

 = \si_{i-1}^{-1} \cdots \si_{1}^{-1} 
   \left( a_{1,2g}^{-1}\cdots  a_{1,1}^{-1}\right)
   \si_{j-1}^{-1} \cdots \si_{2}^{-1} \si_{1}^{-2} \si_2 \cdots \si_{j-1}
   \left( a_{1,1}\cdots a_{1,2g} \right)
   \si_1 \cdots \si_{i-1} \\ \\

 = \si_{i-1}^{-1} \cdots \si_{1}^{-1}\si_{j-1}^{-1} \cdots \si_{2}^{-1}
   \left( a_{1,2g}^{-1}\cdots  a_{1,1}^{-1}\right)  \si_{1}^{-2}
   \left( a_{1,1}\cdots a_{1,2g} \right)
   \si_2 \cdots \si_{j-1} \si_1 \cdots \si_{i-1} \\ \\

 = \left( \si_{i-1}^{-1} \cdots \si_{1}^{-1}
          \si_{j-1}^{-1} \cdots \si_{2}^{-1} \right)
     a_{1,2g}^{-1} A_{1,2g}^{-1} \si_{1}^{-2} A_{1,2g} a_{1,2g} 
   \left( \si_2 \cdots \si_{j-1} \si_1 \cdots \si_{i-1} \right) \\ \\

\mbox{(since $\left( A_{1,2g} A_{2,2g}\right)$ commutes with $\si_1$)} \\ \\

 = \left( \si_{i-1}^{-1} \cdots \si_{1}^{-1}
          \si_{j-1}^{-1} \cdots \si_{2}^{-1} \right)
    a_{1,2g}^{-1} A_{2,2g} \si_{1}^{-2} A_{2,2g}^{-1}  a_{1,2g}
   \left( \si_2 \cdots \si_{j-1} \si_1 \cdots \si_{i-1} \right) \\ \\

\mbox{(by (11))} \\ \\

 = \left( \si_{i-1}^{-1} \cdots \si_{1}^{-1}
          \si_{j-1}^{-1} \cdots \si_{2}^{-1} \right)
 A_{2,2g}  a_{1,2g}^{-1} \si_{1}^{-2}  a_{1,2g} A_{2,2g}^{-1}
   \left( \si_2 \cdots \si_{j-1} \si_1 \cdots \si_{i-1} \right) \\ \\

\mbox{(since, $\left( a_{1,2g} a_{2,2g}\right)$ commutes with 
$\si_1$)} \\ \\

 =  \left( \si_{i-1}^{-1} \cdots \si_{1}^{-1}
          \si_{j-1}^{-1} \cdots \si_{2}^{-1} \right)
 A_{2,2g} a_{2,2g} \si_{1}^{-2} a_{2,2g}^{-1} A_{2,2g}^{-1}
   \left( \si_2 \cdots \si_{j-1} \si_1 \cdots \si_{i-1} \right) 
\end{array}
$$
$$
\begin{array}{l}

 =  \left( \si_{i-1}^{-1} \cdots \si_{1}^{-1} \right)
   A_{j,2g} a_{j,2g} \left( \si_{j-1}^{-1} \cdots \si_{2}^{-1}
  \si_{1}^{-2} \si_2 \cdots \si_{j-1} \right)
  a_{j,2g}^{-1} A_{j,2g}^{-1}  \left( \si_1 \cdots \si_{i-1} \right) \\ \\

\mbox{(by (10))} 
\end{array}
$$
$$
\begin{array}{l}
 =  A_{j,2g} a_{j,2g} \left( \si_{i-1}^{-1} \cdots \si_{1}^{-1} \right)
 \left( \si_{j-1}^{-1} \cdots \si_{2}^{-1}
  \si_{1}^{-2} \si_2 \cdots \si_{j-1} \right)
 \left( \si_1 \cdots \si_{i-1} \right)
 a_{j,2g}^{-1} A_{j,2g}^{-1} \\ \\

\mbox{(by (9))} \\ \\

 =  a_{j,1} \cdots a_{j,2g} 
\left( \si_{j-1}^{-1} \cdots \si_{i+1}^{-1}
  \si_{i}^{-2} \si_{i+1} \cdots \si_{j-1} \right)
  a_{j,2g}^{-1} \cdots a_{j,1}^{-1}.
\end{array}
$$

\vspace{.3cm}
 And this clearly yields (PR8):
$$
\begin{array}{l}
 \left( \prod_{i=1}^{j-1}{a_{i,2g}^{-1}\cdots 
        a_{i,1}^{-1} T_{i,j-1} 
        T_{i,j}^{-1} a_{i,1}\cdots a_{i,2g}}\right) a_{j,1}\cdots 
        a_{j,2g}a_{j,1}^{-1}\cdots a_{j,2g}^{-1} \\ \\

 =  \left( \prod_{i=1}^{j-1}{a_{j,1} \cdots a_{j,2g} 
\left( \si_{j-1}^{-1} \cdots \si_{i+1}^{-1}
  \si_{i}^{-2} \si_{i+1} \cdots \si_{j-1} \right)
  a_{j,2g}^{-1} \cdots a_{j,1}^{-1}}\right) a_{j,1}\cdots 
        a_{j,2g}a_{j,1}^{-1}\cdots a_{j,2g}^{-1} \\ \\

 = a_{j,1} \cdots a_{j,2g} 
  \left( \si_{j-1}^{-1} \cdots \si_{2}^{-1}
  \si_{1}^{-2} \si_{2}^{-1} \cdots \si_{j-1}^{-1} \right)
  a_{j,1}^{-1}\cdots a_{j,2g}^{-1} \\ \\

 = \left( \si_j \cdots \si_{n-1}\right) 
  a_{n,1} \cdots a_{n,2g} \left( \si_{n-1}^{-1} \cdots 
  \si_{1}^{-2} \cdots \si_{n-1}^{-1} \right)
  a_{n,1}^{-1}\cdots a_{n,2g}^{-1}
   \left(\si_{n-1}\cdots \si_j \right) \\ \\

\mbox{(by (R9) and (PR1))} \\ \\

   = \left( \si_j \cdots \si_{n-1}\right)  
    \left(\si_{n-1}\cdots \si_j \right) =  T_{j,n}.
\end{array}
$$

\vspace{.3cm}
 We have thus finished with relations of Type $1$. 

\vspace{.3cm}
  Consider now 
those of Type $2$. For each relator in the presentation of $\Si_n$,
we must find the word $U$ mentioned above.

\vspace{.3cm}
  The first relator is $\de_i \de_j \de_i^{-1} \de_j^{-1}$, when 
$|i-j|\geq 2$ which, by (R1), yields in $\overline{B_n(M)}$ the 
relation
$$
   \si_i \si_j \si_i^{-1} \si_j^{-1}=1 \quad \quad  \quad (|i-j|\geq 2).
$$
Clearly, $U$ is the trivial word.

\vspace{.3cm}
  The second relator, $\de_i \de_{i+1}\de_{i}\de_{i+1}^{-1}
\de_{i}^{-1}\de_{i+1}^{-1}$ gives, by (R2),
$$
  \si_i \si_{i+1}\si_{i}\si_{i+1}^{-1}
\si_{i}^{-1}\si_{i+1}^{-1}=1   \quad \quad (i=1,\cdots,n-2),
$$ 
so in this case $U$ is also the trivial word.

\vspace{.3cm}
  Finally, by the third relator $\de_i^{2}$, we obtain, using (R9),
$$
  \si_i^2 = T_{i,i+1}  \quad  \quad \quad (i=1,\cdots, n-1),
$$
hence $U=T_{i,i+1}$.

\vspace{.3cm}
  So we have obtained the relations in $\overline{B_n(M)}$ mapped
by $\psi$ to the relations of Type~2. 

\vspace{.3cm}
  We finish the proof of
Theorem~\ref{presBnM} obtaining the relations of Type $3$. 
They are very easy to deduce, using (10), (R1), (R2), (R7), (R8) and
(R9). They are the following:
\begin{itemize}
 \item[] $ \si_i a_{j,r} \si_i^{-1} = a_{j,r}    
	\quad \quad \quad \quad \quad (j\neq i, i+1)$, 

 \item[] $ \si_i a_{i,r} \si_i^{-1} = a_{i+1,r} T_{i,i+1}^{-1}
   \quad \quad   \mbox{if $r$ is even,}$

 \item[] $\si_i a_{i,r} \si_i^{-1} = T_{i,i+1} a_{i+1,r} 
  \quad \quad    \mbox{if $r$ is odd,}$

 \item[] $\si_i a_{i+1,r} \si_i^{-1} = T_{i,i+1} a_{i,r} \quad \quad 
    \mbox{if $r$ is even,}$

 \item[] $\si_i a_{i+1,r} \si_i^{-1} = a_{i,r} T_{i,i+1}^{-1} 
  \quad \quad \mbox{if $r$ is odd,}$

 \item[] $ \si_i T_{j,k} \si_i^{-1} = T_{j,k}  \quad \quad 
    \quad \quad \quad (i\neq j-1, j, k),$

 \item[]  $ \si_i T_{i+1,k} \si_i^{-1} = T_{i,k} T_{i,i+1}^{-1},$

 \item[] $  \si_i T_{i,k} \si_i^{-1} = T_{i,i+1} T_{i+1,k},$

 \item[]  $\si_i T_{j,i} \si_i^{-1} = T_{j,i-1} T_{j,i}^{-1} T_{j,i+1}$.
\end{itemize}

\vspace{.6cm}
\section{The braid groups of a non-orientable surface}\label{secnonori}
This section is devoted to prove Theorem~\ref{presBnN}, using the same method
as before. Thus, let $M$ be a closed, non-orientable surface of genus 
$g\geq 2$.

\vspace{.3cm}
\noindent {\bf Step 1.} Denote by $\overline{PB_n(M)}$ the group defined by
the following presentation.

\vspace{.3cm}
\underline{\bf Presentation 3} 

\begin{itemize}

 \item Generators: $\left\{a_{i,r}; \; \; 1\leq i \leq n, \; 1\leq r \leq g\right\}
    \cup  \left\{ T_{j,k}; \;\; 1\leq j < k \leq n \right\}$.

 \item Relations:

 \begin{itemize}

  \item[(Pr1)] $ a_{n,1}^{2}\cdots a_{n,g}^{2}=
        \prod_{i=1}^{n-1}{T_{i,n-1}^{-1}T_{i,n}}. $

  \item[(Pr2)] $  a_{i,r} A_{j,s}= A_{j,s} a_{i,r}$ \hspace{\stretch{1}}
  $ (1\leq i < j \leq n; \: 1\leq r,s \leq g; \: r\neq s ).$  

  \item[(Pr3)] $\left( a_{i,1}^{2}\cdots a_{i,r-1}^{2} a_{i,r}\right) A_{j,r}
     \left( a_{i,r}^{-1}a_{i,r-1}^{-2}\cdots a_{i,1}^{-2} \right) 
     A_{j,r}^{-1}= T_{i,j} T_{i,j-1}^{-1}$ 

   \hspace{\stretch{1}}     $   (1\leq i < j \leq n; \:  1\leq r\leq g). $

  \item[(Pr4)] $ T_{i,j} T_{k,l}=T_{k,l} T_{i,j} $   \hspace{\stretch{1}}
     $  (1 \leq i<j<k<l\leq n \; \mbox{ or } \;1\leq i<k<l \leq j \leq  n ).$

  \item[(Pr5)] $T_{k,l} T_{i,j} T_{k,l}^{-1}=T_{i,k-1}T_{i,k}^{-1}T_{i,j}
         T_{i,l}^{-1} T_{i,k}T_{i,k-1}^{-1}T_{i,l}$ \hspace{\stretch{1}}
      $  (1\leq i<k \leq j<l\leq n).  $

  \item[(Pr6)] $a_{i,r} T_{j,k}=T_{j,k} a_{i,r}$ \hspace{\stretch{1}} 
   $ (1 \leq i<j<k\leq n \; \mbox{ or } \; 1\leq j<k<i\leq n), \: 
    (1\leq r \leq g).$

  \item[(Pr7)] $ a_{i,r} \left( a_{j,g}^{-2}\cdots a_{j,1}^{-2} T_{j,k} \right)
          = \left( a_{j,g}^{-2}\cdots a_{j,1}^{-2} T_{j,k} \right) a_{i,r}$ 
      \hspace{\stretch{1}} $ (1\leq j<i\leq k\leq n).$

  \item[(Pr8)] $ T_{j,n}=a_{j,1}^{2}\cdots a_{j,g}^{2}
        \left( \prod_{i=1}^{j-1}{ T_{j-i,j}^{-1} T_{j-i,j-1}} \right).$

 \end{itemize}

\end{itemize}
Where 
$$
    A_{j,r}=a_{j,1}^{2}\cdots a_{j,r-1}^2 a_{j,r}^{-1}
            a_{j,r-1}^{-2}\cdots a_{j,1}^{-2}.
$$

\vspace{.3cm}
  We shall need, as in the orientable case, another presentation of 
$\overline{PB_n(M)}$, which is the following one.

\vspace{.3cm}
\underline{\bf Presentation 4} 

\begin{itemize}

 \item Generators: $\left\{A_{i,r}; \; \; 1\leq i \leq n, \; 1\leq r \leq g\right\}
    \cup  \left\{ T_{j,k}; \;\; 1\leq j < k \leq n \right\}$.

 \item Relations: the same as in Presentation $3$, where 
$$
   a_{i,r}=A_{i,1}^2\cdots A_{i,r-1}^2 A_{i,r}^{-1} A_{i,r-1}^{-2}
           A_{i,1}^{-2}. 
$$
\end{itemize}

\begin{figure}[ht]
\centerline{\input{nonair.pstex_t} \hspace{1cm} \input{nonTij.pstex_t}}
\caption{The generators of  $PB_n(M)$.}
\label{nonair}
\end {figure}
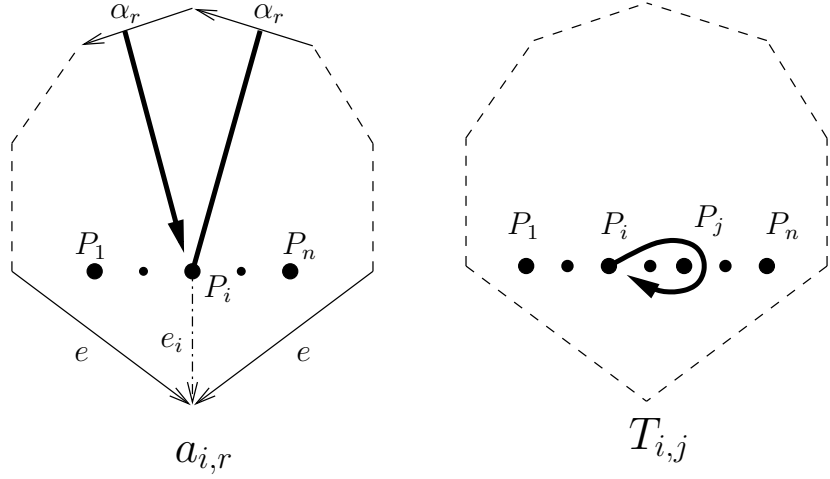

\vspace{.3cm}
 It is clear that Presentation 3 and Presentation 4 are equivalent, in
the same way as they were Presentation 1 and Presentation 2. We must 
now define the homomorphism
$$
 \overline{PB_n(M)} \; \stackrel{\varphi}{\longrightarrow} \; PB_n(M),
$$ 
  by giving the image of the generators. They will be similar 
to those of the orientable surface. For all $i$ and $j$ such that
$1\leq i\leq j \leq n$, the braid $T_{i,j}$ will be the same as in
Section~\ref{secori}.
For all $i$, $r$, such that  $1\leq i\leq n$ and $1\leq r \leq g$,
the braid $a_{i,r}$ will represent the $i$-th string passing through the
$r$-th wall, in the way of Figure~\ref{nonair}. We define as well the path
$e_i$ $(i=1,\ldots,n)$, which goes from $P_i$ to the final point of $e$,
as in Figure~\ref{nonair}. 

\begin{figure}[ht]
\centerline{\input{nonAjr.pstex_t}}
\caption{The braid $A_{j,r}$ in the $P_i$-polygon ($i<j$).}
\label{nonAjr}
\end {figure}
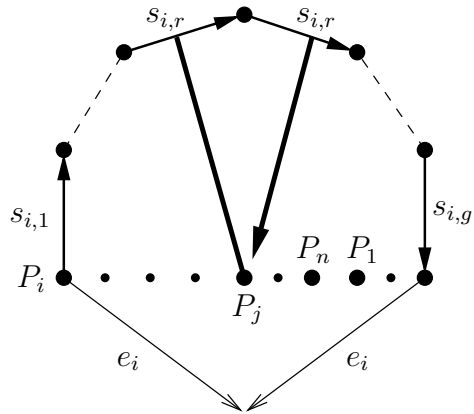

\vspace{.3cm}
  Given $i\in\{1,\ldots,n \}$, denote by $s_{i,r}$ the $i$-th string
of $a_{i,r}$. We can proceed as we did for the $P_1$-polygon in 
Section~\ref{statements} to get the $P_i$-polygon: Cut along the paths 
$e_i$ and $s_{i,1},\ldots, s_{i,g}$, and
glue along $e$ and $\al_1,\ldots, \al_g$. The resulting  $P_i$-polygon is 
labeled by the paths
$$
   s_{i,1}, s_{i,1}, s_{i,2}, s_{i,2},\ldots, s_{i,g}, s_{i,g}, 
  e_i, e_i^{-1},
$$
reading clockwise. Now we can repeat the process of 
Section~\ref{statements} to see that for $1\leq i < j$, the braid 
$$
    A_{j,r}=a_{j,1}^{2}\cdots a_{j,r-1}^2 a_{j,r}^{-1}
            a_{j,r-1}^{-2}\cdots a_{j,1}^{-2}
$$
can be represented in the $P_i$-polygon in the way of 
Figure~\ref{nonAjr}.

\vspace{.3cm}
 The remainder of Step $1$, that is to show that $\varphi$ is a well defined
homomorphism, is analogous to the orientable case. That is, 
Relations (Pr4), (Pr5) and (Pr6) are obvious; Relations (Pr1), (Pr2)
and (Pr3) are analogous to Relations (r3), (r4) and (r5) of 
Theorem~\ref{presBnN}; and we can easily check Relations (Pr7) and (Pr8)
in the $P_j$-polygon.

\vspace{.3cm}
\noindent {\bf Step 2.} This step parallels, up to evident substitutions,
the corresponding one in Section~\ref{secori},
showing the following theorem:

\vspace{.3cm}
\begin{theorem}\label{presPBnN}
  If $M$ is a closed, non-orientable surface of genus $g\geq 2$, then 
$PB_n(M)$ admits Presentation~$3$ (and also Presentation~$4$) as 
presentation. 
\end{theorem}

\vspace{.3cm}
\noindent {\bf Step 3.} 
  Denote by $\overline{B_n(M)}$ the group defined by the presentation
of Theorem~\ref{presBnN}. Call $a_{1,r}$ the elements $a_{r}$ for  
$r=1,\ldots,g$, and then add the generators

\begin{itemize}
 \item[$\relbar$] $ a_{i,r}  \quad \quad \quad i=2,\ldots,n\;; 
     \quad r=1,\ldots,g, $

 \item[$\relbar$] $ T_{j,k}  \quad \quad \quad  1\leq j < k \leq n,  $
\end{itemize}
and the relations
\begin{itemize}

 \item[(r7)] $a_{j+1,r}=\si_j^{-1} a_{j,r}\si_j $   \hspace{\stretch{1}}
      ($1\leq j \leq n-1; \; 1\leq r \leq g$).

 \item[(r8)]  $ T_{j,k}= \si_{j} \si_{j+1}\cdots \si_{k-2} 
          \si_{k-1}^2  \si_{k-2} \cdots \si_j $ \hspace{\stretch{1}}
      $(1 \leq j < k \leq n)$.

\end{itemize}

\vspace{.3cm}
  This provides an equivalent presentation of $\overline{B_n(M)}$,
and the naturally defined function 
$$
  \psi: \overline{B_n(M)} \longrightarrow B_n(M),
$$ 
which is easily proved to be a well defined homomorphism.

\vspace{.3cm}
  Now it remains to apply Lemma~\ref{baselem} to the exact sequence (1),
and then to find relations in $\overline{B_n(M)}$ mapping by
$\psi$ to those of Types 1, 2 and 3, as we did in Section~\ref{secori}.
For those of Type 1 corresponding to (Pr1)-(Pr6), we can use almost the 
same calculations that in the previous section. 

\vspace{.3cm}
  The relation mapping to (Pr7) is obtained as follows:
$$
\begin{array}{l}
   a_{1,r} \left( a_{j,g}^{-2}\cdots a_{j,1}^{-2} T_{j,k} \right) \\ \\

 = a_{i,r} \left( \si_{j-1}^{-1} \cdots \si_1^{-1} 
   a_{1,g}^{-2}\cdots a_{1,1}^{-2} \si_1 \cdots \si_{k-1} \right)
   \si_{k-1} \cdots \si_j \\ \\

\mbox{(by (r3))} \\ \\

 = a_{i,r} \left( \si_{j-1}^{-1} \cdots \si_1^{-1} \right)
     \left( \si_{1}^{-1} \cdots \si_{n-1}^{-1} \right)
     \left( \si_{n-1}^{-1} \cdots \si_{k}^{-1} \right)
     \left( \si_{k-1} \cdots \si_{j} \right) \\ \\

  =  \left( \si_{j-1}^{-1} \cdots \si_1^{-1} \right)
\left( \si_{1}^{-1} \cdots \si_{i-2}^{-1} \right) a_{i,r} \si_{i-1}^{-1}
 \left( \si_{i}^{-1} \cdots \si_{n-1}^{-1} \right)
     \left( \si_{n-1}^{-1} \cdots \si_{k}^{-1} \right)
     \left( \si_{k-1} \cdots \si_{j} \right) \\ \\

=  \left( \si_{j-1}^{-1} \cdots \si_1^{-1} \right)
\left( \si_{1}^{-1} \cdots \si_{i-2}^{-1} \right) 
 \si_{i-1}^{-1} a_{i-1,r}
 \left( \si_{i}^{-1} \cdots \si_{n-1}^{-1} \right)
     \left( \si_{n-1}^{-1} \cdots \si_{k}^{-1} \right)
     \left( \si_{k-1} \cdots \si_{j} \right)  \\ \\

 = \left( \si_{j-1}^{-1} \cdots \si_1^{-1} \right)
     \left( \si_{1}^{-1} \cdots \si_{n-1}^{-1} \right)
\left( \si_{n-1}^{-1} \cdots \si_{k}^{-1} \right)
     \left( \si_{k-1} \cdots \si_{i} \right) a_{i-1,r}
     \si_{i-1}\si_{i-2} \cdots \si_{j} \\ \\

 = \left( \si_{j-1}^{-1} \cdots \si_1^{-1} \right)
     \left( \si_{1}^{-1} \cdots \si_{n-1}^{-1} \right)
     \left( \si_{n-1}^{-1} \cdots \si_{k}^{-1} \right)
     \left( \si_{k-1} \cdots \si_{j} \right) a_{i,r} \\ \\

  = \left( a_{j,g}^{-2}\cdots a_{j,1}^{-2} T_{j,k} \right) a_{i,r},
\end{array}
$$
and the relation mapping to (Pr8), comes from the following calculation:
$$
\hspace{-8cm}
\begin{array}{l}
     a_{j,1}^{2}\cdots a_{j,g}^{2}
        \left( \prod_{i=1}^{j-1}{ T_{j-i,j}^{-1} T_{j-i,j-1}} \right) \\ \\

\mbox{(by (9))} \\ \\

 =  a_{j,1}^{2}\cdots a_{j,g}^{2} \left( \si_{j-1}^{-1} \cdots \si_1^{-2}
  \cdots \si_{j-1}^{-1} \right) 
\end{array}
$$
$$
\hspace{-8cm}
\begin{array}{l}
 
 = \si_{j-1}^{-1} \cdots \si_1^{-1} a_{1,1}^{2}\cdots a_{1,g}^{2}
     \si_1^{-1}  \cdots \si_{j-1}^{-1} \\ \\

\mbox{(by (r3))} \\ \\ 

 = \si_{j} \cdots \si_{n-1}^{2}
  \cdots \si_{j} = T_{j,n}.
\end{array}
$$ 

\vspace{.3cm}
 Finally, the relations mapping by $\psi$ to those of Type 2, are
identical to those for the orientable surfaces,  and relations of Type 3 are
equally easy to deduce. Therefore, we have finished 
the proof of Theorem~\ref{presBnN}.

\vspace{.6cm}
\section{The word problem}\label{word}

   In this section we explain an algorithm to solve the word problem in
the braid group of a surface, using our new presentations. We shall
only explain the orientable case, remarking that the same method can be
used in the non-orientable one.

\vspace{.3cm}
  Let $\om$ be a word over the generators of $B_n(M)$, that is, over $\si_1,
\ldots,\si_{n-1}, a_{1},\ldots,a_{2g}$ and their inverses. The 
algorithm we propose shall give as output a word 
$$
   \om'= \om_1\cdots \om_n s
$$
equivalent to $w$, where $\om_i$ will be a word over $\{ a_{i,1},\ldots,
a_{i,2g},T_{i,i+1},\ldots,T_{i,n-1}\}$, and $s$ will be a word over
$\{\si_1,\ldots,\si_{n-1} \}$ representing the permutation which $\om$ 
induces on the strings. Moreover, we will show that this expression 
is unique, thus $\om=1$ if and only if $\om'$ is the trivial word. 
This algorithm is analogous to the classical braid combing in the braid group
of the disc.

\vspace{.3cm}
  First we need some previous results. Consider the homomorphism $f$ in the
exact sequence (1); it sends $\om$ to its corresponding permutation. 
Now for any element of $\Si_n$, we can take a normal form as a word
over $\{\de_1,\ldots,\de_{n-1} \}$. For instance, we can use the
normal forms in \cite{dlharpe}, where any element of $\Si_n$ is written
as a product 
$$
   t_{1,k_1} t_{2,k_2} \cdots t_{n-1,k_{n-1}}, 
$$
where $t_{m,0}=1$ and  $t_{m,k}= \de_m \de_{m-1} \cdots \de_{m-k+1}$.
 If we replace in this normal form
$\de_i$ by $\si_i$ for $i=1,\ldots, n-1$, we obtain a map
$g: \Si_n \rightarrow W$, where $W$ is the set of words over 
$\{\si_1,\ldots,\si_{n-1} \}$ and their inverses.

 \vspace{.3cm}
  Consider then the composition $\ep= g \circ f$:
$$
   \ep : \; B_n(M) \stackrel{f}{\longrightarrow} \Si_n 
       \stackrel{g}{\longrightarrow} W.
$$
 This map sends any braid to a braid word inducing the same permutation
on the strings. Moreover, the image of $\ep$ is finite, since so is 
$\Si_n$. 

\vspace{.3cm}
  Now in order to apply the algorithm, we need to make a ``dictionary'',
in the following way: for all braid words $p$ in the image of 
$\ep$, consider all braids of the form
$$
   p \; a_{r}^{\pm 1} \; p^{-1}, \quad \quad \quad \quad 
   p \; \si_i^{\pm 1} \; \ep(p\si_i)^{-1}.
$$
Clearly, there is only a finite number of them, and they are all pure
braids. It is not difficult to write these braids as words over 
$\{ a_{i,r}, T_{j,k}\}$ using the relations of the given presentation 
of $B_n(M)$. These are the first words in our dictionary.

\vspace{.3cm}
  Now for $j=1,\ldots,n$, we define the following sets:
$$
   W_j=\{a_{i,r}^{\pm 1}; \;\;i=1,\ldots, j,\; r=1,\ldots,2g\}\cup 
\{T_{i,k}^{\pm 1}; \;\; i=1,\ldots,j, \; k=i+1,\ldots, n-1\},
$$ 
$$
  V_j=\{A_{j,r}^{\pm 1}; \;\; r=1,\ldots, 2g\}\cup 
\{T_{j,k}^{\pm 1} ; \;\; k=j+1,\ldots,n-1\}.
$$ 
For each $x\in W_i$ and each $y\in V_j$, $i<j$, we want to add to 
our dictionary an expression of the form 
$$
   y \; x \; y^{-1} = Z,
$$
where $Z$ is a word over $W_i$. If $y$ is a positive letter, this 
expression is just a relation of Type 3. It may happen that in $Z$ there
is a letter of the form $T_{l,n}^{\pm1}$ ($l\leq i$), but we can replace it
by a word over $W_i$ using (PR8). If $y$ is a negative letter, we can deduce 
the above expression in the same way that we did for relations of Type 3. 
So in any case, we can add all of them to our dictionary.

\vspace{.3cm}
  We still need one more result: Denote by $S_{n,r}$ the $n$-th 
string of $A_{n,r}$. Since $s_{n,1},\ldots,s_{n,2g}$ generates 
$\pi_1(M, P_n)$, Lemma~\ref{xlem} clearly implies that $\{S_{n,1},\ldots,
S_{n,2g}\}$ is another set of generators. Moreover, applying the formulae
of Lemma~\ref{xlem}, one has
$$
   s_{n,1}^{-1}s_{n,2}^{-1}\cdots s_{n,2g}^{-1}
        s_{n,1}s_{n,2}\cdots s_{n,2g} = 
  \left( S_{n,2g}^{-1}S_{n,2g-1} S_{n,2g-2}^{-1}\cdots S_{n,1} \right)
\left( S_{n,2g} S_{n,2g-1}^{-1} S_{n,2g-2}\cdots S_{n,1}^{-1}\right). 
$$
Hence we obtain:
$$
   \pi_1(M,P_n)=  \left< \{S_{n,1},\ldots, S_{n,2g}\}; \; 
 \left( S_{n,2g}^{-1}S_{n,2g-1} \cdots S_{n,2}^{-1} S_{n,1} \right) 
 \left(  S_{n,2g} S_{n,2g-1}^{-1} \cdots S_{n,2} S_{n,1}^{-1}
 \right)=1 \right>.
$$

\vspace{.3cm} 
  We are now ready to start with the algorithm. Thus, let $\om$ be
a word over the generators of $B_n(M)$. Define the word $s = \ep(\om)$. 
Since the normal forms in $\Si_n$ are unique, so is $s$. We obtain 
a word $\overline{\om}= \om s^{-1} \in PB_n(M)$ such that $\om=\overline{\om} s$.

\vspace{.3cm}
  Next we want to write $\overline{\om}$ as a word over 
$\{ a_{i,r}, T_{j,k}\}$ (where $i$, $r$, $j$ and $k$ take all possible 
values).  Suppose that $\overline{\om}$ has length $m$, that is,
$\overline{\om}=x_1\cdots x_m$ where each $x_i$ is a 
generator of $B_n(M)$ or its inverse. For all
$i=1,\ldots,m$ define $\overline{\om}_i=x_1\cdots x_i$. 
Since $\overline{\om}\in PB_n(M)$, then 
$\ep(\overline{\om}_m)= \ep(\overline{\om})=1$, so one has:
$$
  \overline{\om}= x_1\cdots x_m 
   =\left( 1 \; x_1 \;\ep(\overline{\om}_1)^{-1}\right)
  \left(\ep(\overline{\om}_1) \; x_2 \;\ep(\overline{\om}_2)^{-1}\right)
  \cdots 
  \left(\ep(\overline{\om}_{m-1})\; x_m \; \ep(\overline{\om}_m)^{-1}\right).
$$
But all factors on the right hand side of the equation are included in
our dictionary, so we can use it to write all of them, and thus
$\overline{\om}$, as a word over the generators of $PB_n(M)$.

\vspace{.3cm}
  The next step is to replace in $\overline{\om}$ all the letters of the 
form $a_{n,r}^{\pm 1}$ using the 
formula in Presentation $2$,
$$
   a_{n,r}^{(-1)^{r+1}}=
  \left( A_{n,1} A_{n,2}^{-1} A_{n,3}\cdots A_{n,r-1}^{\pm 1}\right)
  \left( A_{n,r+1}^{\mp 1}  \cdots A_{n,2g-1}^{-1} A_{n,2g}\right),
$$
and all the letters of the form $T_{j,n}^{\pm1}$, using (PR8).
In this way we obtain $\overline{\om}$ written as a word over
$W_{n-1}\cup V_n$. We use again the dictionary  to ``move'' to the 
right hand side of $\overline{\om}$ all the letters in $V_n$. We 
will obtain $\overline{\om}= X\; Y$, where  $X$ is a word over 
$W_{n-1}$ and $Y$ is a word over $V_n$.

\vspace{.3cm}
 Consider now the following exact sequence, coming also from the
Fadell-Neuwirth fibration (see~\cite{birman}).
$$
     1\longrightarrow PB_{n-1}(M\backslash \{ P_n\}) 
      \stackrel{u}{\longrightarrow} PB_n(M) 
    \stackrel{v}{\longrightarrow} \pi_1(M,P_n) \longrightarrow 1,
$$
where for all $\Ga=(\ga_1,\ldots,\ga_n)\in PB_n(M)$, $\; v(\Ga)=\ga_n$. 
Note that 
$v(\overline{\om})= Y \in \pi_1(M)$. Now in $\pi_1(M)$ we could apply Dehn's
algorithm (see~\cite{lyndonschupp}) to obtain a normal form of $Y$. 
At each step of Dehn's algorithm, a sub-word of $Y$ would be replaced by
a shorter one, using the relation 
$$
\left( S_{n,2g}^{-1}S_{n,2g-1} S_{n,2g-2}^{-1}\cdots S_{n,1}\right)
\left(  S_{n,2g} S_{n,2g-1}^{-1} S_{n,2g-2}\cdots S_{n,1}^{-1}
 \right)=1.
$$
Instead of this, we will do a similar process in $PB_n(M)$: each time that 
Dehn's algorithm replaces a sub-word of $Y$ in $\pi_1(M)$, we replace the 
corresponding sub-word in $\overline{\om}=XY \in PB_n(M)$ using
$$
 \left( A_{n,2g}^{-1}A_{n,2g-1} A_{n,2g-2}^{-1}\cdots A_{n,1}\right)
\left(  A_{n,2g} A_{n,2g-1}^{-1} A_{n,2g-2}\cdots A_{n,1}^{-1}
 \right) = \prod_{i=1}^{n-1}{T_{i,n-1}^{-1}T_{i,n}},
$$
which is a relation equivalent to (PR1); then we remove the $T_{i,n}^{\pm 1}$
using (PR8) and we move again the letters in $V_n$ to the right hand side
of our word.

\vspace{.3cm}
  At the end of this process, we will obtain $\overline{\om}= X_{n-1}\; 
\om_n$, where $\om_n$ is the normal form of $v(\overline{\om})$ in 
$\pi_1(M)$, so it is unique, and $X_{n-1}$ is a word over 
$W_{n-1}$.

\vspace{.3cm}
  The algorithm will end in $n-1$ steps: 
At each step, we have a word $X_{m}$ over  $W_{m}$, we replace the letters
of the form $a_{m,r}^{\pm 1}$ by words over $V_m$, and then we 
move all the letters of $V_m$ to the right hand side, using the dictionary.
Then we remove all the sub-words of the form $x x^{-1}$ or $x^{-1} x$,
and we obtain $X_{m}=X_{m-1}\om_{m}$, where $X_{m-1}$ is a word over 
$W_{m-1}$ and $\om_{m}$ is a reduced word over $V_m$. If we prove that 
the word $\om_{m}$ is unique, we will have the unique factorization
$\om = \om_1\cdots \om_n s$ as the output of our algorithm.

\vspace{.3cm}
  Define $M_{n-m}=M\backslash \{ P_{m+1},\ldots,P_n\}$ for any 
$m=1,\ldots, n-1$. In~\cite{birman} we can find the following
exact sequence, analogous to the previous one.
$$
     1\longrightarrow PB_{m-1}(M_{n-m+1}) 
      \stackrel{f}{\longrightarrow} PB_m(M_{n-m}) 
    \stackrel{g}{\longrightarrow} \pi_1(M_{n-m}) 
\longrightarrow 1.
$$
We only need to notice that 
$X_{m}\in PB_m(M_{n-m})$, and 
$g(X_{m})=\om_{m}$. Now since $\pi_1(M_{n-m})$
is a free group with free system of generators $\{a_{m,r}; \; 
1\leq r \leq 2g \}\cup \{T_{m,j}; \; m+1 \leq j\leq n-1\}$, and since 
$\om_{m}$ is a reduced word, then it is unique,
as we wanted to show.

\vspace{.4cm}
\begin{tabular}{ll}
\noindent J. GONZ\'ALEZ-MENESES \\
Universit\'e de Bourgogne \hspace{4truecm} & Departamento de \'Algebra \\
Laboratoire de Topologie & Facultad de Matem\'aticas\\
UMR 5584 du CNRS & Universidad de Sevilla\\
B. P. 47870  & C/ Tarfia, s/n\\
21078 - Dijon Cedex (France)&  41012 - Sevilla (Spain)\\
{\em jmeneses@u-bourgogne.fr}& {\em meneses@algebra.us.es}
\end{tabular}

\end{document}

%% file: define.tex
\newcommand\al{\alpha}

\newcommand\ga{\gamma}
\newcommand\de{\delta}
\newcommand\ep{\varepsilon}

\newcommand\si{\sigma}

\newcommand\om{\omega}

\newcommand\Ga{\Gamma}
\newcommand\De{\Delta}

\newcommand\Si{\Sigma}

\DeclareMathAlphabet{\frak}{U}{euf}{m}{n}
\SetMathAlphabet{\frak}{bold}{U}{euf}{b}{n}
\DeclareMathAlphabet{\Bbb}{U}{msb}{m}{n}

  \newcounter{numero}[section]
  \newcounter{aux1}
  \newcounter{aux2}
  \renewcommand{\thenumero}{\thesection.\arabic{numero}}
  \newenvironment{theorem}{
  \dimen255=\parindent \parindent=0in
  \refstepcounter{numero}{\vspace{.3cm}\bf Theorem \thenumero .}\
  \it}{\rm \parindent=\dimen255\par}

  \newenvironment{lemma}{
  \dimen255=\parindent \parindent=0in
  \refstepcounter{numero} {\vspace{.3cm}\bf Lemma \thenumero .}\
  \it}{\rm \parindent=\dimen255\par}

  \newenvironment{proof}{
  {\vspace{.3cm}\noindent\sc{Proof:}}\
  }{$\fbox{}$ \par}

%% file: polygon.pstex_t
\begin{picture}(0,0)%
\includegraphics{polygon.pstex}%
\end{picture}%
\setlength{\unitlength}{3947sp}%
\begingroup\makeatletter\ifx\SetFigFont\undefined
\def\x#1#2#3#4#5#6#7\relax{\def\x{#1#2#3#4#5#6}}%
\expandafter\x\fmtname xxxxxx\relax \def\y{splain}%
\ifx\x\y   
\gdef\SetFigFont#1#2#3{%
  \ifnum #1<17\tiny\else \ifnum #1<20\small\else
  \ifnum #1<24\normalsize\else \ifnum #1<29\large\else
  \ifnum #1<34\Large\else \ifnum #1<41\LARGE\else
     \huge\fi\fi\fi\fi\fi\fi
  \csname #3\endcsname}%
\else
\gdef\SetFigFont#1#2#3{\begingroup
  \count@#1\relax \ifnum 25<\count@\count@25\fi
  \def\x{\endgroup\@setsize\SetFigFont{#2pt}}%
  \expandafter\x
    \csname \romannumeral\the\count@ pt\expandafter\endcsname
    \csname @\romannumeral\the\count@ pt\endcsname
  \csname #3\endcsname}%
\fi
\fi\endgroup
\begin{picture}(2700,2383)(496,-2063)
\put(3196,-511){\makebox(0,0)[lb]{\smash{\SetFigFont{12}{14.4}{rm}$\al_1$}}}
\put(496,-421){\makebox(0,0)[lb]{\smash{\SetFigFont{12}{14.4}{rm}$\al_{2g}$}}}
\put(2746,-1951){\makebox(0,0)[lb]{\smash{\SetFigFont{12}{14.4}{rm}$\al_{2g-1}$}}}
\put(631,-1366){\makebox(0,0)[lb]{\smash{\SetFigFont{12}{14.4}{rm}$\al_1$}}}
\put(1081,-1996){\makebox(0,0)[lb]{\smash{\SetFigFont{12}{14.4}{rm}$\al_2$}}}
\put(766,164){\makebox(0,0)[lb]{\smash{\SetFigFont{12}{14.4}{rm}$\al_{2g-1}$}}}
\put(3196,-1321){\makebox(0,0)[lb]{\smash{\SetFigFont{12}{14.4}{rm}$\al_{2g}$}}}
\put(2746,119){\makebox(0,0)[lb]{\smash{\SetFigFont{12}{14.4}{rm}$\al_2$}}}
\end{picture}

%% file: viewpoint1.pstex_t
\begin{picture}(0,0)%
\includegraphics{viewpoint1.pstex}%
\end{picture}%
\setlength{\unitlength}{3947sp}%
\begingroup\makeatletter\ifx\SetFigFont\undefined
\def\x#1#2#3#4#5#6#7\relax{\def\x{#1#2#3#4#5#6}}%
\expandafter\x\fmtname xxxxxx\relax \def\y{splain}%
\ifx\x\y   
\gdef\SetFigFont#1#2#3{%
  \ifnum #1<17\tiny\else \ifnum #1<20\small\else
  \ifnum #1<24\normalsize\else \ifnum #1<29\large\else
  \ifnum #1<34\Large\else \ifnum #1<41\LARGE\else
     \huge\fi\fi\fi\fi\fi\fi
  \csname #3\endcsname}%
\else
\gdef\SetFigFont#1#2#3{\begingroup
  \count@#1\relax \ifnum 25<\count@\count@25\fi
  \def\x{\endgroup\@setsize\SetFigFont{#2pt}}%
  \expandafter\x
    \csname \romannumeral\the\count@ pt\expandafter\endcsname
    \csname @\romannumeral\the\count@ pt\endcsname
  \csname #3\endcsname}%
\fi
\fi\endgroup
\begin{picture}(2667,2082)(271,-1906)
\put(1801,-1906){\makebox(0,0)[lb]{\smash{\SetFigFont{12}{14.4}{rm}$L$}}}
\put(271,-691){\makebox(0,0)[lb]{\smash{\SetFigFont{12}{14.4}{rm}$I$}}}
\end{picture}

%% file: viewpoint2.pstex_t
\begin{picture}(0,0)%
\includegraphics{viewpoint2.pstex}%
\end{picture}%
\setlength{\unitlength}{3947sp}%
\begingroup\makeatletter\ifx\SetFigFont\undefined
\def\x#1#2#3#4#5#6#7\relax{\def\x{#1#2#3#4#5#6}}%
\expandafter\x\fmtname xxxxxx\relax \def\y{splain}%
\ifx\x\y   
\gdef\SetFigFont#1#2#3{%
  \ifnum #1<17\tiny\else \ifnum #1<20\small\else
  \ifnum #1<24\normalsize\else \ifnum #1<29\large\else
  \ifnum #1<34\Large\else \ifnum #1<41\LARGE\else
     \huge\fi\fi\fi\fi\fi\fi
  \csname #3\endcsname}%
\else
\gdef\SetFigFont#1#2#3{\begingroup
  \count@#1\relax \ifnum 25<\count@\count@25\fi
  \def\x{\endgroup\@setsize\SetFigFont{#2pt}}%
  \expandafter\x
    \csname \romannumeral\the\count@ pt\expandafter\endcsname
    \csname @\romannumeral\the\count@ pt\endcsname
  \csname #3\endcsname}%
\fi
\fi\endgroup
\begin{picture}(2364,2364)(635,-2136)
\end{picture}

%% file: a1r.pstex_t
\begin{picture}(0,0)%
\includegraphics{a1r.pstex}%
\end{picture}%
\setlength{\unitlength}{3947sp}%
\begingroup\makeatletter\ifx\SetFigFont\undefined
\def\x#1#2#3#4#5#6#7\relax{\def\x{#1#2#3#4#5#6}}%
\expandafter\x\fmtname xxxxxx\relax \def\y{splain}%
\ifx\x\y   
\gdef\SetFigFont#1#2#3{%
  \ifnum #1<17\tiny\else \ifnum #1<20\small\else
  \ifnum #1<24\normalsize\else \ifnum #1<29\large\else
  \ifnum #1<34\Large\else \ifnum #1<41\LARGE\else
     \huge\fi\fi\fi\fi\fi\fi
  \csname #3\endcsname}%
\else
\gdef\SetFigFont#1#2#3{\begingroup
  \count@#1\relax \ifnum 25<\count@\count@25\fi
  \def\x{\endgroup\@setsize\SetFigFont{#2pt}}%
  \expandafter\x
    \csname \romannumeral\the\count@ pt\expandafter\endcsname
    \csname @\romannumeral\the\count@ pt\endcsname
  \csname #3\endcsname}%
\fi
\fi\endgroup
\begin{picture}(2409,2985)(619,-2491)
\put(1621,-2491){\makebox(0,0)[lb]{\smash{\SetFigFont{17}{20.4}{rm}$a_{2k+1}$}}}
\put(2296,-2131){\makebox(0,0)[lb]{\smash{\SetFigFont{12}{14.4}{rm}$\al_{2k+1}$}}}
\put(1171,299){\makebox(0,0)[lb]{\smash{\SetFigFont{12}{14.4}{rm}$\al_{2k+1}$}}}
\put(2611,-691){\makebox(0,0)[lb]{\smash{\SetFigFont{12}{14.4}{rm}$P_{n}$}}}
\put(1621,-691){\makebox(0,0)[lb]{\smash{\SetFigFont{12}{14.4}{rm}$P_{1}$}}}
\end{picture}

%% file: a1s.pstex_t
\begin{picture}(0,0)%
\includegraphics{a1s.pstex}%
\end{picture}%
\setlength{\unitlength}{3947sp}%
\begingroup\makeatletter\ifx\SetFigFont\undefined
\def\x#1#2#3#4#5#6#7\relax{\def\x{#1#2#3#4#5#6}}%
\expandafter\x\fmtname xxxxxx\relax \def\y{splain}%
\ifx\x\y   
\gdef\SetFigFont#1#2#3{%
  \ifnum #1<17\tiny\else \ifnum #1<20\small\else
  \ifnum #1<24\normalsize\else \ifnum #1<29\large\else
  \ifnum #1<34\Large\else \ifnum #1<41\LARGE\else
     \huge\fi\fi\fi\fi\fi\fi
  \csname #3\endcsname}%
\else
\gdef\SetFigFont#1#2#3{\begingroup
  \count@#1\relax \ifnum 25<\count@\count@25\fi
  \def\x{\endgroup\@setsize\SetFigFont{#2pt}}%
  \expandafter\x
    \csname \romannumeral\the\count@ pt\expandafter\endcsname
    \csname @\romannumeral\the\count@ pt\endcsname
  \csname #3\endcsname}%
\fi
\fi\endgroup
\begin{picture}(2409,2985)(619,-2491)
\put(1351,-691){\makebox(0,0)[lb]{\smash{\SetFigFont{12}{14.4}{rm}$P_{1}$}}}
\put(2386,299){\makebox(0,0)[lb]{\smash{\SetFigFont{12}{14.4}{rm}$\al_{2k}$}}}
\put(1261,-2086){\makebox(0,0)[lb]{\smash{\SetFigFont{12}{14.4}{rm}$\al_{2k}$}}}
\put(2611,-691){\makebox(0,0)[lb]{\smash{\SetFigFont{12}{14.4}{rm}$P_{n}$}}}
\put(1711,-2491){\makebox(0,0)[lb]{\smash{\SetFigFont{17}{20.4}{rm}$a_{2k}$}}}
\end{picture}

%% file: sii.pstex_t
\begin{picture}(0,0)%
\includegraphics{sii.pstex}%
\end{picture}%
\setlength{\unitlength}{3947sp}%
\begingroup\makeatletter\ifx\SetFigFont\undefined
\def\x#1#2#3#4#5#6#7\relax{\def\x{#1#2#3#4#5#6}}%
\expandafter\x\fmtname xxxxxx\relax \def\y{splain}%
\ifx\x\y   
\gdef\SetFigFont#1#2#3{%
  \ifnum #1<17\tiny\else \ifnum #1<20\small\else
  \ifnum #1<24\normalsize\else \ifnum #1<29\large\else
  \ifnum #1<34\Large\else \ifnum #1<41\LARGE\else
     \huge\fi\fi\fi\fi\fi\fi
  \csname #3\endcsname}%
\else
\gdef\SetFigFont#1#2#3{\begingroup
  \count@#1\relax \ifnum 25<\count@\count@25\fi
  \def\x{\endgroup\@setsize\SetFigFont{#2pt}}%
  \expandafter\x
    \csname \romannumeral\the\count@ pt\expandafter\endcsname
    \csname @\romannumeral\the\count@ pt\endcsname
  \csname #3\endcsname}%
\fi
\fi\endgroup
\begin{picture}(2409,2712)(619,-2491)
\put(1801,-2491){\makebox(0,0)[lb]{\smash{\SetFigFont{17}{20.4}{rm}$\si_{i}$}}}
\put(2701,-691){\makebox(0,0)[lb]{\smash{\SetFigFont{12}{14.4}{rm}$P_{n}$}}}
\put(2206,-691){\makebox(0,0)[lb]{\smash{\SetFigFont{12}{14.4}{rm}$P_{i+1}$}}}
\put(946,-691){\makebox(0,0)[lb]{\smash{\SetFigFont{12}{14.4}{rm}$P_{1}$}}}
\put(1486,-691){\makebox(0,0)[lb]{\smash{\SetFigFont{12}{14.4}{rm}$P_{i}$}}}
\end{picture}

%% file: puzzleP11.pstex_t
\begin{picture}(0,0)%
\includegraphics{puzzleP11.pstex}%
\end{picture}%
\setlength{\unitlength}{3947sp}%
\begingroup\makeatletter\ifx\SetFigFont\undefined
\def\x#1#2#3#4#5#6#7\relax{\def\x{#1#2#3#4#5#6}}%
\expandafter\x\fmtname xxxxxx\relax \def\y{splain}%
\ifx\x\y   
\gdef\SetFigFont#1#2#3{%
  \ifnum #1<17\tiny\else \ifnum #1<20\small\else
  \ifnum #1<24\normalsize\else \ifnum #1<29\large\else
  \ifnum #1<34\Large\else \ifnum #1<41\LARGE\else
     \huge\fi\fi\fi\fi\fi\fi
  \csname #3\endcsname}%
\else
\gdef\SetFigFont#1#2#3{\begingroup
  \count@#1\relax \ifnum 25<\count@\count@25\fi
  \def\x{\endgroup\@setsize\SetFigFont{#2pt}}%
  \expandafter\x
    \csname \romannumeral\the\count@ pt\expandafter\endcsname
    \csname @\romannumeral\the\count@ pt\endcsname
  \csname #3\endcsname}%
\fi
\fi\endgroup
\begin{picture}(2520,2637)(541,-2143)
\put(1441,-916){\makebox(0,0)[lb]{\smash{\SetFigFont{12}{14.4}{rm}$P_{1}$}}}
\put(2386,-2086){\makebox(0,0)[lb]{\smash{\SetFigFont{12}{14.4}{rm}$\al_{3}$}}}
\put(1306,299){\makebox(0,0)[lb]{\smash{\SetFigFont{12}{14.4}{rm}$\al_{3}$}}}
\put(2341,299){\makebox(0,0)[lb]{\smash{\SetFigFont{12}{14.4}{rm}$\al_{2}$}}}
\put(1261,-2086){\makebox(0,0)[lb]{\smash{\SetFigFont{12}{14.4}{rm}$\al_{2}$}}}
\put(541,-421){\makebox(0,0)[lb]{\smash{\SetFigFont{12}{14.4}{rm}$\al_4$}}}
\put(541,-1456){\makebox(0,0)[lb]{\smash{\SetFigFont{12}{14.4}{rm}$\al_1$}}}
\put(3061,-1456){\makebox(0,0)[lb]{\smash{\SetFigFont{12}{14.4}{rm}$\al_4$}}}
\put(3061,-466){\makebox(0,0)[lb]{\smash{\SetFigFont{12}{14.4}{rm}$\al_1$}}}
\put(2746,-916){\makebox(0,0)[lb]{\smash{\SetFigFont{12}{14.4}{rm}$P_{n}$}}}
\put(1036,-511){\makebox(0,0)[lb]{\smash{\SetFigFont{12}{14.4}{rm}$s_{4}$}}}
\put(2206,-286){\makebox(0,0)[lb]{\smash{\SetFigFont{12}{14.4}{rm}$s_{2}$}}}
\put(1396,-286){\makebox(0,0)[lb]{\smash{\SetFigFont{12}{14.4}{rm}$s_{3}$}}}
\put(2431,-511){\makebox(0,0)[lb]{\smash{\SetFigFont{12}{14.4}{rm}$s_{1}$}}}
\end{picture}

%% file: puzzleP12.pstex_t
\begin{picture}(0,0)%
\includegraphics{puzzleP12.pstex}%
\end{picture}%
\setlength{\unitlength}{3947sp}%
\begingroup\makeatletter\ifx\SetFigFont\undefined
\def\x#1#2#3#4#5#6#7\relax{\def\x{#1#2#3#4#5#6}}%
\expandafter\x\fmtname xxxxxx\relax \def\y{splain}%
\ifx\x\y   
\gdef\SetFigFont#1#2#3{%
  \ifnum #1<17\tiny\else \ifnum #1<20\small\else
  \ifnum #1<24\normalsize\else \ifnum #1<29\large\else
  \ifnum #1<34\Large\else \ifnum #1<41\LARGE\else
     \huge\fi\fi\fi\fi\fi\fi
  \csname #3\endcsname}%
\else
\gdef\SetFigFont#1#2#3{\begingroup
  \count@#1\relax \ifnum 25<\count@\count@25\fi
  \def\x{\endgroup\@setsize\SetFigFont{#2pt}}%
  \expandafter\x
    \csname \romannumeral\the\count@ pt\expandafter\endcsname
    \csname @\romannumeral\the\count@ pt\endcsname
  \csname #3\endcsname}%
\fi
\fi\endgroup
\begin{picture}(2880,3030)(361,-2356)
\put(1846,-2356){\makebox(0,0)[lb]{\smash{\SetFigFont{12}{14.4}{rm}$P_{1}$}}}
\put(2386,-2086){\makebox(0,0)[lb]{\smash{\SetFigFont{12}{14.4}{rm}$s_{2}$}}}
\put(2341,299){\makebox(0,0)[lb]{\smash{\SetFigFont{12}{14.4}{rm}$s_{3}$}}}
\put(3061,-1456){\makebox(0,0)[lb]{\smash{\SetFigFont{12}{14.4}{rm}$s_{1}$}}}
\put(3061,-466){\makebox(0,0)[lb]{\smash{\SetFigFont{12}{14.4}{rm}$s_{4}$}}}
\put(1036,-511){\makebox(0,0)[lb]{\smash{\SetFigFont{12}{14.4}{rm}$\al_{1}$}}}
\put(2431,-511){\makebox(0,0)[lb]{\smash{\SetFigFont{12}{14.4}{rm}$\al_{4}$}}}
\put(3241,-961){\makebox(0,0)[lb]{\smash{\SetFigFont{12}{14.4}{rm}$P_{1}$}}}
\put(1441,-916){\makebox(0,0)[lb]{\smash{\SetFigFont{12}{14.4}{rm}$P_{n}$}}}
\put(361,-916){\makebox(0,0)[lb]{\smash{\SetFigFont{12}{14.4}{rm}$P_{1}$}}}
\put(1441,-286){\makebox(0,0)[lb]{\smash{\SetFigFont{12}{14.4}{rm}$\al_{2}$}}}
\put(2251,-286){\makebox(0,0)[lb]{\smash{\SetFigFont{12}{14.4}{rm}$\al_{3}$}}}
\put(1351,299){\makebox(0,0)[lb]{\smash{\SetFigFont{12}{14.4}{rm}$s_{2}$}}}
\put(586,-421){\makebox(0,0)[lb]{\smash{\SetFigFont{12}{14.4}{rm}$s_{1}$}}}
\put(586,-1456){\makebox(0,0)[lb]{\smash{\SetFigFont{12}{14.4}{rm}$s_{4}$}}}
\put(1306,-2086){\makebox(0,0)[lb]{\smash{\SetFigFont{12}{14.4}{rm}$s_{3}$}}}
\put(811,119){\makebox(0,0)[lb]{\smash{\SetFigFont{12}{14.4}{rm}$P_{1}$}}}
\put(1846,479){\makebox(0,0)[lb]{\smash{\SetFigFont{12}{14.4}{rm}$P_{1}$}}}
\put(2836,119){\makebox(0,0)[lb]{\smash{\SetFigFont{12}{14.4}{rm}$P_{1}$}}}
\put(766,-1906){\makebox(0,0)[lb]{\smash{\SetFigFont{12}{14.4}{rm}$P_{1}$}}}
\put(2836,-1906){\makebox(0,0)[lb]{\smash{\SetFigFont{12}{14.4}{rm}$P_{1}$}}}
\end{picture}

%% file: B3.pstex_t
\begin{picture}(0,0)%
\includegraphics{B3.pstex}%
\end{picture}%
\setlength{\unitlength}{3947sp}%
\begingroup\makeatletter\ifx\SetFigFont\undefined
\def\x#1#2#3#4#5#6#7\relax{\def\x{#1#2#3#4#5#6}}%
\expandafter\x\fmtname xxxxxx\relax \def\y{splain}%
\ifx\x\y   
\gdef\SetFigFont#1#2#3{%
  \ifnum #1<17\tiny\else \ifnum #1<20\small\else
  \ifnum #1<24\normalsize\else \ifnum #1<29\large\else
  \ifnum #1<34\Large\else \ifnum #1<41\LARGE\else
     \huge\fi\fi\fi\fi\fi\fi
  \csname #3\endcsname}%
\else
\gdef\SetFigFont#1#2#3{\begingroup
  \count@#1\relax \ifnum 25<\count@\count@25\fi
  \def\x{\endgroup\@setsize\SetFigFont{#2pt}}%
  \expandafter\x
    \csname \romannumeral\the\count@ pt\expandafter\endcsname
    \csname @\romannumeral\the\count@ pt\endcsname
  \csname #3\endcsname}%
\fi
\fi\endgroup
\begin{picture}(2880,2384)(361,-2153)
\put(1711,-916){\makebox(0,0)[lb]{\smash{\SetFigFont{12}{14.4}{rm}$P_{n}$}}}
\put(3061,-1456){\makebox(0,0)[lb]{\smash{\SetFigFont{12}{14.4}{rm}$s_{1}$}}}
\put(3061,-466){\makebox(0,0)[lb]{\smash{\SetFigFont{12}{14.4}{rm}$s_{2g}$}}}
\put(3241,-961){\makebox(0,0)[lb]{\smash{\SetFigFont{12}{14.4}{rm}$P_{1}$}}}
\put(361,-916){\makebox(0,0)[lb]{\smash{\SetFigFont{12}{14.4}{rm}$P_{1}$}}}
\put(541,-1456){\makebox(0,0)[lb]{\smash{\SetFigFont{12}{14.4}{rm}$s_{2g}$}}}
\put(586,-421){\makebox(0,0)[lb]{\smash{\SetFigFont{12}{14.4}{rm}$s_{1}$}}}
\end{picture}

%% file: A2r2.pstex_t
\begin{picture}(0,0)%
\includegraphics{A2r2.pstex}%
\end{picture}%
\setlength{\unitlength}{3947sp}%
\begingroup\makeatletter\ifx\SetFigFont\undefined
\def\x#1#2#3#4#5#6#7\relax{\def\x{#1#2#3#4#5#6}}%
\expandafter\x\fmtname xxxxxx\relax \def\y{splain}%
\ifx\x\y   
\gdef\SetFigFont#1#2#3{%
  \ifnum #1<17\tiny\else \ifnum #1<20\small\else
  \ifnum #1<24\normalsize\else \ifnum #1<29\large\else
  \ifnum #1<34\Large\else \ifnum #1<41\LARGE\else
     \huge\fi\fi\fi\fi\fi\fi
  \csname #3\endcsname}%
\else
\gdef\SetFigFont#1#2#3{\begingroup
  \count@#1\relax \ifnum 25<\count@\count@25\fi
  \def\x{\endgroup\@setsize\SetFigFont{#2pt}}%
  \expandafter\x
    \csname \romannumeral\the\count@ pt\expandafter\endcsname
    \csname @\romannumeral\the\count@ pt\endcsname
  \csname #3\endcsname}%
\fi
\fi\endgroup
\begin{picture}(3240,3210)(136,-2401)
\put(1981,-781){\makebox(0,0)[lb]{\smash{\SetFigFont{12}{14.4}{rm}$P_2$}}}
\put(3376,-691){\makebox(0,0)[lb]{\smash{\SetFigFont{12}{14.4}{rm}$\al_{1}$}}}
\put(3241,-241){\makebox(0,0)[lb]{\smash{\SetFigFont{12}{14.4}{rm}$\al_{2}$}}}
\put(2206,-2356){\makebox(0,0)[lb]{\smash{\SetFigFont{12}{14.4}{rm}$\al_{r+1}$}}}
\put(1216,-2356){\makebox(0,0)[lb]{\smash{\SetFigFont{12}{14.4}{rm}$\al_{r-1}$}}}
\put(1171,569){\makebox(0,0)[lb]{\smash{\SetFigFont{12}{14.4}{rm}$\al_{r+1}$}}}
\put(2206,569){\makebox(0,0)[lb]{\smash{\SetFigFont{12}{14.4}{rm}$\al_{r-1}$}}}
\put(136,-691){\makebox(0,0)[lb]{\smash{\SetFigFont{12}{14.4}{rm}$\al_{2g}$}}}
\put(3016,-1906){\makebox(0,0)[lb]{\smash{\SetFigFont{12}{14.4}{rm}$\al_{2g-2}$}}}
\put(3061, 74){\makebox(0,0)[lb]{\smash{\SetFigFont{12}{14.4}{rm}$\al_{3}$}}}
\put(1621,-781){\makebox(0,0)[lb]{\smash{\SetFigFont{12}{14.4}{rm}$P_1$}}}
\put(226,-1096){\makebox(0,0)[lb]{\smash{\SetFigFont{12}{14.4}{rm}$\al_{1}$}}}
\put(316,-1501){\makebox(0,0)[lb]{\smash{\SetFigFont{12}{14.4}{rm}$\al_{2}$}}}
\put(541,-1861){\makebox(0,0)[lb]{\smash{\SetFigFont{12}{14.4}{rm}$\al_{3}$}}}
\put(3376,-1096){\makebox(0,0)[lb]{\smash{\SetFigFont{12}{14.4}{rm}$\al_{2g}$}}}
\put(3241,-1546){\makebox(0,0)[lb]{\smash{\SetFigFont{12}{14.4}{rm}$\al_{2g-1}$}}}
\put(136,-241){\makebox(0,0)[lb]{\smash{\SetFigFont{12}{14.4}{rm}$\al_{2g-1}$}}}
\put(361,119){\makebox(0,0)[lb]{\smash{\SetFigFont{12}{14.4}{rm}$\al_{2g-2}$}}}
\put(1846,614){\makebox(0,0)[lb]{\smash{\SetFigFont{12}{14.4}{rm}$\al_{r}$}}}
\put(1801,-2401){\makebox(0,0)[lb]{\smash{\SetFigFont{12}{14.4}{rm}$\al_{r}$}}}
\end{picture}

%% file: A2r1.pstex_t
\begin{picture}(0,0)%
\includegraphics{A2r1.pstex}%
\end{picture}%
\setlength{\unitlength}{3947sp}%
\begingroup\makeatletter\ifx\SetFigFont\undefined
\def\x#1#2#3#4#5#6#7\relax{\def\x{#1#2#3#4#5#6}}%
\expandafter\x\fmtname xxxxxx\relax \def\y{splain}%
\ifx\x\y   
\gdef\SetFigFont#1#2#3{%
  \ifnum #1<17\tiny\else \ifnum #1<20\small\else
  \ifnum #1<24\normalsize\else \ifnum #1<29\large\else
  \ifnum #1<34\Large\else \ifnum #1<41\LARGE\else
     \huge\fi\fi\fi\fi\fi\fi
  \csname #3\endcsname}%
\else
\gdef\SetFigFont#1#2#3{\begingroup
  \count@#1\relax \ifnum 25<\count@\count@25\fi
  \def\x{\endgroup\@setsize\SetFigFont{#2pt}}%
  \expandafter\x
    \csname \romannumeral\the\count@ pt\expandafter\endcsname
    \csname @\romannumeral\the\count@ pt\endcsname
  \csname #3\endcsname}%
\fi
\fi\endgroup
\begin{picture}(3285,3210)(91,-2401)
\put(676,-916){\makebox(0,0)[lb]{\smash{\SetFigFont{12}{14.4}{rm}$P_2$}}}
\put(1801,614){\makebox(0,0)[lb]{\smash{\SetFigFont{12}{14.4}{rm}$s_{r}$}}}
\put(3376,-691){\makebox(0,0)[lb]{\smash{\SetFigFont{12}{14.4}{rm}$s_{2g}$}}}
\put(3241,-241){\makebox(0,0)[lb]{\smash{\SetFigFont{12}{14.4}{rm}$s_{2g-1}$}}}
\put(3376,-1141){\makebox(0,0)[lb]{\smash{\SetFigFont{12}{14.4}{rm}$s_{1}$}}}
\put(3241,-1591){\makebox(0,0)[lb]{\smash{\SetFigFont{12}{14.4}{rm}$s_{2}$}}}
\put(271,-286){\makebox(0,0)[lb]{\smash{\SetFigFont{12}{14.4}{rm}$s_{2}$}}}
\put(136,-691){\makebox(0,0)[lb]{\smash{\SetFigFont{12}{14.4}{rm}$s_{1}$}}}
\put( 91,-1096){\makebox(0,0)[lb]{\smash{\SetFigFont{12}{14.4}{rm}$s_{2g}$}}}
\put( 91,-1501){\makebox(0,0)[lb]{\smash{\SetFigFont{12}{14.4}{rm}$s_{2g-1}$}}}
\put(2206,-2356){\makebox(0,0)[lb]{\smash{\SetFigFont{12}{14.4}{rm}$s_{r-1}$}}}
\put(1216,-2356){\makebox(0,0)[lb]{\smash{\SetFigFont{12}{14.4}{rm}$s_{r+1}$}}}
\put(1756,-2401){\makebox(0,0)[lb]{\smash{\SetFigFont{12}{14.4}{rm}$s_{r}$}}}
\put(226,-916){\makebox(0,0)[lb]{\smash{\SetFigFont{12}{14.4}{rm}$P_1$}}}
\put(1171,569){\makebox(0,0)[lb]{\smash{\SetFigFont{12}{14.4}{rm}$s_{r-1}$}}}
\put(2206,569){\makebox(0,0)[lb]{\smash{\SetFigFont{12}{14.4}{rm}$s_{r+1}$}}}
\end{picture}

%% file: B51.pstex_t
\begin{picture}(0,0)%
\includegraphics{B51.pstex}%
\end{picture}%
\setlength{\unitlength}{3947sp}%
\begingroup\makeatletter\ifx\SetFigFont\undefined
\def\x#1#2#3#4#5#6#7\relax{\def\x{#1#2#3#4#5#6}}%
\expandafter\x\fmtname xxxxxx\relax \def\y{splain}%
\ifx\x\y   
\gdef\SetFigFont#1#2#3{%
  \ifnum #1<17\tiny\else \ifnum #1<20\small\else
  \ifnum #1<24\normalsize\else \ifnum #1<29\large\else
  \ifnum #1<34\Large\else \ifnum #1<41\LARGE\else
     \huge\fi\fi\fi\fi\fi\fi
  \csname #3\endcsname}%
\else
\gdef\SetFigFont#1#2#3{\begingroup
  \count@#1\relax \ifnum 25<\count@\count@25\fi
  \def\x{\endgroup\@setsize\SetFigFont{#2pt}}%
  \expandafter\x
    \csname \romannumeral\the\count@ pt\expandafter\endcsname
    \csname @\romannumeral\the\count@ pt\endcsname
  \csname #3\endcsname}%
\fi
\fi\endgroup
\begin{picture}(2632,2805)(361,-2221)
\put(361,-916){\makebox(0,0)[lb]{\smash{\SetFigFont{12}{14.4}{rm}$P_1$}}}
\put(1936,-1006){\makebox(0,0)[lb]{\smash{\SetFigFont{12}{14.4}{rm}$P_2$}}}
\put(1801,389){\makebox(0,0)[lb]{\smash{\SetFigFont{12}{14.4}{rm}$s_{r}$}}}
\put(1801,-2221){\makebox(0,0)[lb]{\smash{\SetFigFont{12}{14.4}{rm}$s_{r}$}}}
\put(541,-556){\makebox(0,0)[lb]{\smash{\SetFigFont{12}{14.4}{rm}$s_{1}$}}}
\end{picture}

%% file: B52.pstex_t
\begin{picture}(0,0)%
\includegraphics{B52.pstex}%
\end{picture}%
\setlength{\unitlength}{3947sp}%
\begingroup\makeatletter\ifx\SetFigFont\undefined
\def\x#1#2#3#4#5#6#7\relax{\def\x{#1#2#3#4#5#6}}%
\expandafter\x\fmtname xxxxxx\relax \def\y{splain}%
\ifx\x\y   
\gdef\SetFigFont#1#2#3{%
  \ifnum #1<17\tiny\else \ifnum #1<20\small\else
  \ifnum #1<24\normalsize\else \ifnum #1<29\large\else
  \ifnum #1<34\Large\else \ifnum #1<41\LARGE\else
     \huge\fi\fi\fi\fi\fi\fi
  \csname #3\endcsname}%
\else
\gdef\SetFigFont#1#2#3{\begingroup
  \count@#1\relax \ifnum 25<\count@\count@25\fi
  \def\x{\endgroup\@setsize\SetFigFont{#2pt}}%
  \expandafter\x
    \csname \romannumeral\the\count@ pt\expandafter\endcsname
    \csname @\romannumeral\the\count@ pt\endcsname
  \csname #3\endcsname}%
\fi
\fi\endgroup
\begin{picture}(2632,2805)(361,-2221)
\put(541,-556){\makebox(0,0)[lb]{\smash{\SetFigFont{12}{14.4}{rm}$s_{1}$}}}
\put(1936,-1006){\makebox(0,0)[lb]{\smash{\SetFigFont{12}{14.4}{rm}$P_2$}}}
\put(361,-916){\makebox(0,0)[lb]{\smash{\SetFigFont{12}{14.4}{rm}$P_1$}}}
\put(1801,389){\makebox(0,0)[lb]{\smash{\SetFigFont{12}{14.4}{rm}$s_{r}$}}}
\put(1801,-2221){\makebox(0,0)[lb]{\smash{\SetFigFont{12}{14.4}{rm}$s_{r}$}}}
\end{picture}

%% file: B53.pstex_t
\begin{picture}(0,0)%
\includegraphics{B53.pstex}%
\end{picture}%
\setlength{\unitlength}{3947sp}%
\begingroup\makeatletter\ifx\SetFigFont\undefined
\def\x#1#2#3#4#5#6#7\relax{\def\x{#1#2#3#4#5#6}}%
\expandafter\x\fmtname xxxxxx\relax \def\y{splain}%
\ifx\x\y   
\gdef\SetFigFont#1#2#3{%
  \ifnum #1<17\tiny\else \ifnum #1<20\small\else
  \ifnum #1<24\normalsize\else \ifnum #1<29\large\else
  \ifnum #1<34\Large\else \ifnum #1<41\LARGE\else
     \huge\fi\fi\fi\fi\fi\fi
  \csname #3\endcsname}%
\else
\gdef\SetFigFont#1#2#3{\begingroup
  \count@#1\relax \ifnum 25<\count@\count@25\fi
  \def\x{\endgroup\@setsize\SetFigFont{#2pt}}%
  \expandafter\x
    \csname \romannumeral\the\count@ pt\expandafter\endcsname
    \csname @\romannumeral\the\count@ pt\endcsname
  \csname #3\endcsname}%
\fi
\fi\endgroup
\begin{picture}(2632,2805)(361,-2221)
\put(1801,-2221){\makebox(0,0)[lb]{\smash{\SetFigFont{12}{14.4}{rm}$s_{r}$}}}
\put(361,-916){\makebox(0,0)[lb]{\smash{\SetFigFont{12}{14.4}{rm}$P_1$}}}
\put(2026,-1051){\makebox(0,0)[lb]{\smash{\SetFigFont{12}{14.4}{rm}$P_2$}}}
\put(541,-556){\makebox(0,0)[lb]{\smash{\SetFigFont{12}{14.4}{rm}$s_{1}$}}}
\put(1801,389){\makebox(0,0)[lb]{\smash{\SetFigFont{12}{14.4}{rm}$s_{r}$}}}
\end{picture}

%% file: nonpolygon.pstex_t
\begin{picture}(0,0)%
\includegraphics{nonpolygon.pstex}%
\end{picture}%
\setlength{\unitlength}{3947sp}%
\begingroup\makeatletter\ifx\SetFigFont\undefined
\def\x#1#2#3#4#5#6#7\relax{\def\x{#1#2#3#4#5#6}}%
\expandafter\x\fmtname xxxxxx\relax \def\y{splain}%
\ifx\x\y   
\gdef\SetFigFont#1#2#3{%
  \ifnum #1<17\tiny\else \ifnum #1<20\small\else
  \ifnum #1<24\normalsize\else \ifnum #1<29\large\else
  \ifnum #1<34\Large\else \ifnum #1<41\LARGE\else
     \huge\fi\fi\fi\fi\fi\fi
  \csname #3\endcsname}%
\else
\gdef\SetFigFont#1#2#3{\begingroup
  \count@#1\relax \ifnum 25<\count@\count@25\fi
  \def\x{\endgroup\@setsize\SetFigFont{#2pt}}%
  \expandafter\x
    \csname \romannumeral\the\count@ pt\expandafter\endcsname
    \csname @\romannumeral\the\count@ pt\endcsname
  \csname #3\endcsname}%
\fi
\fi\endgroup
\begin{picture}(2835,2465)(271,-2144)
\put(2881,-106){\makebox(0,0)[lb]{\smash{\SetFigFont{12}{14.4}{rm}$\al_{2}$}}}
\put(2881,-1681){\makebox(0,0)[lb]{\smash{\SetFigFont{12}{14.4}{rm}$\al_{1}$}}}
\put(1936,-1501){\makebox(0,0)[lb]{\smash{\SetFigFont{12}{14.4}{rm}$e$}}}
\put(2206,-2086){\makebox(0,0)[lb]{\smash{\SetFigFont{12}{14.4}{rm}$\al_1$}}}
\put(496,-61){\makebox(0,0)[lb]{\smash{\SetFigFont{12}{14.4}{rm}$\al_{g-1}$}}}
\put(271,-871){\makebox(0,0)[lb]{\smash{\SetFigFont{12}{14.4}{rm}$\al_{g-1}$}}}
\put(631,-1591){\makebox(0,0)[lb]{\smash{\SetFigFont{12}{14.4}{rm}$\al_{g}$}}}
\put(1306,-2086){\makebox(0,0)[lb]{\smash{\SetFigFont{12}{14.4}{rm}$\al_{g}$}}}
\put(3106,-871){\makebox(0,0)[lb]{\smash{\SetFigFont{12}{14.4}{rm}$\al_2$}}}
\end{picture}

%% file: nonpolygon2.pstex_t
\begin{picture}(0,0)%
\includegraphics{nonpolygon2.pstex}%
\end{picture}%
\setlength{\unitlength}{3947sp}%
\begingroup\makeatletter\ifx\SetFigFont\undefined
\def\x#1#2#3#4#5#6#7\relax{\def\x{#1#2#3#4#5#6}}%
\expandafter\x\fmtname xxxxxx\relax \def\y{splain}%
\ifx\x\y   
\gdef\SetFigFont#1#2#3{%
  \ifnum #1<17\tiny\else \ifnum #1<20\small\else
  \ifnum #1<24\normalsize\else \ifnum #1<29\large\else
  \ifnum #1<34\Large\else \ifnum #1<41\LARGE\else
     \huge\fi\fi\fi\fi\fi\fi
  \csname #3\endcsname}%
\else
\gdef\SetFigFont#1#2#3{\begingroup
  \count@#1\relax \ifnum 25<\count@\count@25\fi
  \def\x{\endgroup\@setsize\SetFigFont{#2pt}}%
  \expandafter\x
    \csname \romannumeral\the\count@ pt\expandafter\endcsname
    \csname @\romannumeral\the\count@ pt\endcsname
  \csname #3\endcsname}%
\fi
\fi\endgroup
\begin{picture}(2655,2419)(406,-2098)
\put(3061,-916){\makebox(0,0)[lb]{\smash{\SetFigFont{12}{14.4}{rm}$\al_1$}}}
\put(2476,-1771){\makebox(0,0)[lb]{\smash{\SetFigFont{12}{14.4}{rm}$e$}}}
\put(1081,-1771){\makebox(0,0)[lb]{\smash{\SetFigFont{12}{14.4}{rm}$e$}}}
\put(2431,-1096){\makebox(0,0)[lb]{\smash{\SetFigFont{12}{14.4}{rm}$P_n$}}}
\put(1126,-1096){\makebox(0,0)[lb]{\smash{\SetFigFont{12}{14.4}{rm}$P_1$}}}
\put(676,-16){\makebox(0,0)[lb]{\smash{\SetFigFont{12}{14.4}{rm}$\al_{g}$}}}
\put(406,-826){\makebox(0,0)[lb]{\smash{\SetFigFont{12}{14.4}{rm}$\al_{g}$}}}
\put(2881,-106){\makebox(0,0)[lb]{\smash{\SetFigFont{12}{14.4}{rm}$\al_{1}$}}}
\end{picture}

%% file: nona1r.pstex_t
\begin{picture}(0,0)%
\includegraphics{nona1r.pstex}%
\end{picture}%
\setlength{\unitlength}{3947sp}%
\begingroup\makeatletter\ifx\SetFigFont\undefined
\def\x#1#2#3#4#5#6#7\relax{\def\x{#1#2#3#4#5#6}}%
\expandafter\x\fmtname xxxxxx\relax \def\y{splain}%
\ifx\x\y   
\gdef\SetFigFont#1#2#3{%
  \ifnum #1<17\tiny\else \ifnum #1<20\small\else
  \ifnum #1<24\normalsize\else \ifnum #1<29\large\else
  \ifnum #1<34\Large\else \ifnum #1<41\LARGE\else
     \huge\fi\fi\fi\fi\fi\fi
  \csname #3\endcsname}%
\else
\gdef\SetFigFont#1#2#3{\begingroup
  \count@#1\relax \ifnum 25<\count@\count@25\fi
  \def\x{\endgroup\@setsize\SetFigFont{#2pt}}%
  \expandafter\x
    \csname \romannumeral\the\count@ pt\expandafter\endcsname
    \csname @\romannumeral\the\count@ pt\endcsname
  \csname #3\endcsname}%
\fi
\fi\endgroup
\begin{picture}(2185,2985)(663,-2401)
\put(1711,-2401){\makebox(0,0)[lb]{\smash{\SetFigFont{17}{20.4}{rm}$a_r$}}}
\put(2476,-1771){\makebox(0,0)[lb]{\smash{\SetFigFont{12}{14.4}{rm}$e$}}}
\put(1081,-1771){\makebox(0,0)[lb]{\smash{\SetFigFont{12}{14.4}{rm}$e$}}}
\put(1306,389){\makebox(0,0)[lb]{\smash{\SetFigFont{12}{14.4}{rm}$\al_{r}$}}}
\put(2206,389){\makebox(0,0)[lb]{\smash{\SetFigFont{12}{14.4}{rm}$\al_r$}}}
\put(2386,-1096){\makebox(0,0)[lb]{\smash{\SetFigFont{12}{14.4}{rm}$P_n$}}}
\put(1666,-1636){\makebox(0,0)[lb]{\smash{\SetFigFont{12}{14.4}{rm}$e_1$}}}
\put(1216,-1096){\makebox(0,0)[lb]{\smash{\SetFigFont{12}{14.4}{rm}$P_1$}}}
\end{picture}

%% file: nonsii.pstex_t
\begin{picture}(0,0)%
\includegraphics{nonsii.pstex}%
\end{picture}%
\setlength{\unitlength}{3947sp}%
\begingroup\makeatletter\ifx\SetFigFont\undefined
\def\x#1#2#3#4#5#6#7\relax{\def\x{#1#2#3#4#5#6}}%
\expandafter\x\fmtname xxxxxx\relax \def\y{splain}%
\ifx\x\y   
\gdef\SetFigFont#1#2#3{%
  \ifnum #1<17\tiny\else \ifnum #1<20\small\else
  \ifnum #1<24\normalsize\else \ifnum #1<29\large\else
  \ifnum #1<34\Large\else \ifnum #1<41\LARGE\else
     \huge\fi\fi\fi\fi\fi\fi
  \csname #3\endcsname}%
\else
\gdef\SetFigFont#1#2#3{\begingroup
  \count@#1\relax \ifnum 25<\count@\count@25\fi
  \def\x{\endgroup\@setsize\SetFigFont{#2pt}}%
  \expandafter\x
    \csname \romannumeral\the\count@ pt\expandafter\endcsname
    \csname @\romannumeral\the\count@ pt\endcsname
  \csname #3\endcsname}%
\fi
\fi\endgroup
\begin{picture}(2655,2722)(406,-2401)
\put(2611,-1006){\makebox(0,0)[lb]{\smash{\SetFigFont{12}{14.4}{rm}$P_{n}$}}}
\put(2116,-1006){\makebox(0,0)[lb]{\smash{\SetFigFont{12}{14.4}{rm}$P_{i+1}$}}}
\put(856,-1006){\makebox(0,0)[lb]{\smash{\SetFigFont{12}{14.4}{rm}$P_{1}$}}}
\put(1396,-1006){\makebox(0,0)[lb]{\smash{\SetFigFont{12}{14.4}{rm}$P_{i}$}}}
\put(1711,-2401){\makebox(0,0)[lb]{\smash{\SetFigFont{17}{20.4}{rm}$\si_i$}}}
\put(2476,-1771){\makebox(0,0)[lb]{\smash{\SetFigFont{12}{14.4}{rm}$e$}}}
\put(1081,-1771){\makebox(0,0)[lb]{\smash{\SetFigFont{12}{14.4}{rm}$e$}}}
\put(676,-16){\makebox(0,0)[lb]{\smash{\SetFigFont{12}{14.4}{rm}$\al_{g}$}}}
\put(406,-826){\makebox(0,0)[lb]{\smash{\SetFigFont{12}{14.4}{rm}$\al_{g}$}}}
\put(2881,-106){\makebox(0,0)[lb]{\smash{\SetFigFont{12}{14.4}{rm}$\al_{1}$}}}
\put(3061,-916){\makebox(0,0)[lb]{\smash{\SetFigFont{12}{14.4}{rm}$\al_1$}}}
\end{picture}

%% file: air.pstex_t
\begin{picture}(0,0)%
\includegraphics{air.pstex}%
\end{picture}%
\setlength{\unitlength}{3947sp}%
\begingroup\makeatletter\ifx\SetFigFont\undefined
\def\x#1#2#3#4#5#6#7\relax{\def\x{#1#2#3#4#5#6}}%
\expandafter\x\fmtname xxxxxx\relax \def\y{splain}%
\ifx\x\y   
\gdef\SetFigFont#1#2#3{%
  \ifnum #1<17\tiny\else \ifnum #1<20\small\else
  \ifnum #1<24\normalsize\else \ifnum #1<29\large\else
  \ifnum #1<34\Large\else \ifnum #1<41\LARGE\else
     \huge\fi\fi\fi\fi\fi\fi
  \csname #3\endcsname}%
\else
\gdef\SetFigFont#1#2#3{\begingroup
  \count@#1\relax \ifnum 25<\count@\count@25\fi
  \def\x{\endgroup\@setsize\SetFigFont{#2pt}}%
  \expandafter\x
    \csname \romannumeral\the\count@ pt\expandafter\endcsname
    \csname @\romannumeral\the\count@ pt\endcsname
  \csname #3\endcsname}%
\fi
\fi\endgroup
\begin{picture}(2409,2985)(619,-2491)
\put(1171,299){\makebox(0,0)[lb]{\smash{\SetFigFont{12}{14.4}{rm}$\al_{2k+1}$}}}
\put(1891,-691){\makebox(0,0)[lb]{\smash{\SetFigFont{12}{14.4}{rm}$P_{i}$}}}
\put(1036,-691){\makebox(0,0)[lb]{\smash{\SetFigFont{12}{14.4}{rm}$P_{1}$}}}
\put(2656,-691){\makebox(0,0)[lb]{\smash{\SetFigFont{12}{14.4}{rm}$P_{n}$}}}
\put(2296,-2131){\makebox(0,0)[lb]{\smash{\SetFigFont{12}{14.4}{rm}$\al_{2k+1}$}}}
\put(1576,-2491){\makebox(0,0)[lb]{\smash{\SetFigFont{17}{20.4}{rm}$a_{i,2k+1}$}}}
\end{picture}

%% file: ais.pstex_t
\begin{picture}(0,0)%
\includegraphics{ais.pstex}%
\end{picture}%
\setlength{\unitlength}{3947sp}%
\begingroup\makeatletter\ifx\SetFigFont\undefined
\def\x#1#2#3#4#5#6#7\relax{\def\x{#1#2#3#4#5#6}}%
\expandafter\x\fmtname xxxxxx\relax \def\y{splain}%
\ifx\x\y   
\gdef\SetFigFont#1#2#3{%
  \ifnum #1<17\tiny\else \ifnum #1<20\small\else
  \ifnum #1<24\normalsize\else \ifnum #1<29\large\else
  \ifnum #1<34\Large\else \ifnum #1<41\LARGE\else
     \huge\fi\fi\fi\fi\fi\fi
  \csname #3\endcsname}%
\else
\gdef\SetFigFont#1#2#3{\begingroup
  \count@#1\relax \ifnum 25<\count@\count@25\fi
  \def\x{\endgroup\@setsize\SetFigFont{#2pt}}%
  \expandafter\x
    \csname \romannumeral\the\count@ pt\expandafter\endcsname
    \csname @\romannumeral\the\count@ pt\endcsname
  \csname #3\endcsname}%
\fi
\fi\endgroup
\begin{picture}(2409,2985)(619,-2491)
\put(1666,-2491){\makebox(0,0)[lb]{\smash{\SetFigFont{17}{20.4}{rm}$a_{i,2k}$}}}
\put(1036,-691){\makebox(0,0)[lb]{\smash{\SetFigFont{12}{14.4}{rm}$P_{1}$}}}
\put(2656,-691){\makebox(0,0)[lb]{\smash{\SetFigFont{12}{14.4}{rm}$P_{n}$}}}
\put(1756,-691){\makebox(0,0)[lb]{\smash{\SetFigFont{12}{14.4}{rm}$P_{i}$}}}
\put(2386,299){\makebox(0,0)[lb]{\smash{\SetFigFont{12}{14.4}{rm}$\al_{2k}$}}}
\put(1261,-2086){\makebox(0,0)[lb]{\smash{\SetFigFont{12}{14.4}{rm}$\al_{2k}$}}}
\end{picture}

%% file: Tij.pstex_t
\begin{picture}(0,0)%
\includegraphics{Tij.pstex}%
\end{picture}%
\setlength{\unitlength}{3947sp}%
\begingroup\makeatletter\ifx\SetFigFont\undefined
\def\x#1#2#3#4#5#6#7\relax{\def\x{#1#2#3#4#5#6}}%
\expandafter\x\fmtname xxxxxx\relax \def\y{splain}%
\ifx\x\y   
\gdef\SetFigFont#1#2#3{%
  \ifnum #1<17\tiny\else \ifnum #1<20\small\else
  \ifnum #1<24\normalsize\else \ifnum #1<29\large\else
  \ifnum #1<34\Large\else \ifnum #1<41\LARGE\else
     \huge\fi\fi\fi\fi\fi\fi
  \csname #3\endcsname}%
\else
\gdef\SetFigFont#1#2#3{\begingroup
  \count@#1\relax \ifnum 25<\count@\count@25\fi
  \def\x{\endgroup\@setsize\SetFigFont{#2pt}}%
  \expandafter\x
    \csname \romannumeral\the\count@ pt\expandafter\endcsname
    \csname @\romannumeral\the\count@ pt\endcsname
  \csname #3\endcsname}%
\fi
\fi\endgroup
\begin{picture}(2409,2703)(619,-2482)
\put(1801,-2401){\makebox(0,0)[lb]{\smash{\SetFigFont{17}{20.4}{rm}$T_{i,j}$}}}
\put(1036,-691){\makebox(0,0)[lb]{\smash{\SetFigFont{12}{14.4}{rm}$P_{1}$}}}
\put(2656,-691){\makebox(0,0)[lb]{\smash{\SetFigFont{12}{14.4}{rm}$P_{n}$}}}
\put(1621,-691){\makebox(0,0)[lb]{\smash{\SetFigFont{12}{14.4}{rm}$P_{i}$}}}
\put(2206,-646){\makebox(0,0)[lb]{\smash{\SetFigFont{12}{14.4}{rm}$P_{j}$}}}
\end{picture}

%% file: Ajr.pstex_t
\begin{picture}(0,0)%
\includegraphics{Ajr.pstex}%
\end{picture}%
\setlength{\unitlength}{3947sp}%
\begingroup\makeatletter\ifx\SetFigFont\undefined
\def\x#1#2#3#4#5#6#7\relax{\def\x{#1#2#3#4#5#6}}%
\expandafter\x\fmtname xxxxxx\relax \def\y{splain}%
\ifx\x\y   
\gdef\SetFigFont#1#2#3{%
  \ifnum #1<17\tiny\else \ifnum #1<20\small\else
  \ifnum #1<24\normalsize\else \ifnum #1<29\large\else
  \ifnum #1<34\Large\else \ifnum #1<41\LARGE\else
     \huge\fi\fi\fi\fi\fi\fi
  \csname #3\endcsname}%
\else
\gdef\SetFigFont#1#2#3{\begingroup
  \count@#1\relax \ifnum 25<\count@\count@25\fi
  \def\x{\endgroup\@setsize\SetFigFont{#2pt}}%
  \expandafter\x
    \csname \romannumeral\the\count@ pt\expandafter\endcsname
    \csname @\romannumeral\the\count@ pt\endcsname
  \csname #3\endcsname}%
\fi
\fi\endgroup
\begin{picture}(2632,2805)(361,-2221)
\put(361,-916){\makebox(0,0)[lb]{\smash{\SetFigFont{12}{14.4}{rm}$P_i$}}}
\put(1936,-1006){\makebox(0,0)[lb]{\smash{\SetFigFont{12}{14.4}{rm}$P_j$}}}
\put(1711,389){\makebox(0,0)[lb]{\smash{\SetFigFont{12}{14.4}{rm}$s_{i,r}$}}}
\put(1711,-2221){\makebox(0,0)[lb]{\smash{\SetFigFont{12}{14.4}{rm}$s_{i,r}$}}}
\end{picture}

%% file: rel81.pstex_t
\begin{picture}(0,0)%
\includegraphics{rel81.pstex}%
\end{picture}%
\setlength{\unitlength}{3947sp}%
\begingroup\makeatletter\ifx\SetFigFont\undefined
\def\x#1#2#3#4#5#6#7\relax{\def\x{#1#2#3#4#5#6}}%
\expandafter\x\fmtname xxxxxx\relax \def\y{splain}%
\ifx\x\y   
\gdef\SetFigFont#1#2#3{%
  \ifnum #1<17\tiny\else \ifnum #1<20\small\else
  \ifnum #1<24\normalsize\else \ifnum #1<29\large\else
  \ifnum #1<34\Large\else \ifnum #1<41\LARGE\else
     \huge\fi\fi\fi\fi\fi\fi
  \csname #3\endcsname}%
\else
\gdef\SetFigFont#1#2#3{\begingroup
  \count@#1\relax \ifnum 25<\count@\count@25\fi
  \def\x{\endgroup\@setsize\SetFigFont{#2pt}}%
  \expandafter\x
    \csname \romannumeral\the\count@ pt\expandafter\endcsname
    \csname @\romannumeral\the\count@ pt\endcsname
  \csname #3\endcsname}%
\fi
\fi\endgroup
\begin{picture}(2835,2294)(226,-2108)
\put(316,-556){\makebox(0,0)[lb]{\smash{\SetFigFont{12}{14.4}{rm}$s_{i,1}$}}}
\put(2251,-1141){\makebox(0,0)[lb]{\smash{\SetFigFont{12}{14.4}{rm}$P_n$}}}
\put(1441,-1141){\makebox(0,0)[lb]{\smash{\SetFigFont{12}{14.4}{rm}$P_j$}}}
\put(226,-961){\makebox(0,0)[lb]{\smash{\SetFigFont{12}{14.4}{rm}$P_i$}}}
\put(3061,-961){\makebox(0,0)[lb]{\smash{\SetFigFont{12}{14.4}{rm}$P_i$}}}
\put(2926,-511){\makebox(0,0)[lb]{\smash{\SetFigFont{12}{14.4}{rm}$s_{i,2g}$}}}
\end{picture}

%% file: rel82.pstex_t
\begin{picture}(0,0)%
\includegraphics{rel82.pstex}%
\end{picture}%
\setlength{\unitlength}{3947sp}%
\begingroup\makeatletter\ifx\SetFigFont\undefined
\def\x#1#2#3#4#5#6#7\relax{\def\x{#1#2#3#4#5#6}}%
\expandafter\x\fmtname xxxxxx\relax \def\y{splain}%
\ifx\x\y   
\gdef\SetFigFont#1#2#3{%
  \ifnum #1<17\tiny\else \ifnum #1<20\small\else
  \ifnum #1<24\normalsize\else \ifnum #1<29\large\else
  \ifnum #1<34\Large\else \ifnum #1<41\LARGE\else
     \huge\fi\fi\fi\fi\fi\fi
  \csname #3\endcsname}%
\else
\gdef\SetFigFont#1#2#3{\begingroup
  \count@#1\relax \ifnum 25<\count@\count@25\fi
  \def\x{\endgroup\@setsize\SetFigFont{#2pt}}%
  \expandafter\x
    \csname \romannumeral\the\count@ pt\expandafter\endcsname
    \csname @\romannumeral\the\count@ pt\endcsname
  \csname #3\endcsname}%
\fi
\fi\endgroup
\begin{picture}(2632,2294)(361,-2108)
\put(1171,-1141){\makebox(0,0)[lb]{\smash{\SetFigFont{12}{14.4}{rm}$P_n$}}}
\put(2026,-1141){\makebox(0,0)[lb]{\smash{\SetFigFont{12}{14.4}{rm}$P_i$}}}
\put(361,-961){\makebox(0,0)[lb]{\smash{\SetFigFont{12}{14.4}{rm}$P_j$}}}
\end{picture}

%% file: nonair.pstex_t
\begin{picture}(0,0)%
\includegraphics{nonair.pstex}%
\end{picture}%
\setlength{\unitlength}{3947sp}%
\begingroup\makeatletter\ifx\SetFigFont\undefined
\def\x#1#2#3#4#5#6#7\relax{\def\x{#1#2#3#4#5#6}}%
\expandafter\x\fmtname xxxxxx\relax \def\y{splain}%
\ifx\x\y   
\gdef\SetFigFont#1#2#3{%
  \ifnum #1<17\tiny\else \ifnum #1<20\small\else
  \ifnum #1<24\normalsize\else \ifnum #1<29\large\else
  \ifnum #1<34\Large\else \ifnum #1<41\LARGE\else
     \huge\fi\fi\fi\fi\fi\fi
  \csname #3\endcsname}%
\else
\gdef\SetFigFont#1#2#3{\begingroup
  \count@#1\relax \ifnum 25<\count@\count@25\fi
  \def\x{\endgroup\@setsize\SetFigFont{#2pt}}%
  \expandafter\x
    \csname \romannumeral\the\count@ pt\expandafter\endcsname
    \csname @\romannumeral\the\count@ pt\endcsname
  \csname #3\endcsname}%
\fi
\fi\endgroup
\begin{picture}(2185,2985)(663,-2401)
\put(1711,-2401){\makebox(0,0)[lb]{\smash{\SetFigFont{17}{20.4}{rm}$a_{i,r}$}}}
\put(2476,-1771){\makebox(0,0)[lb]{\smash{\SetFigFont{12}{14.4}{rm}$e$}}}
\put(1081,-1771){\makebox(0,0)[lb]{\smash{\SetFigFont{12}{14.4}{rm}$e$}}}
\put(1306,389){\makebox(0,0)[lb]{\smash{\SetFigFont{12}{14.4}{rm}$\al_{r}$}}}
\put(2206,389){\makebox(0,0)[lb]{\smash{\SetFigFont{12}{14.4}{rm}$\al_r$}}}
\put(1621,-1681){\makebox(0,0)[lb]{\smash{\SetFigFont{12}{14.4}{rm}$e_i$}}}
\put(1891,-1366){\makebox(0,0)[lb]{\smash{\SetFigFont{12}{14.4}{rm}$P_i$}}}
\put(2386,-1096){\makebox(0,0)[lb]{\smash{\SetFigFont{12}{14.4}{rm}$P_n$}}}
\put(1081,-1096){\makebox(0,0)[lb]{\smash{\SetFigFont{12}{14.4}{rm}$P_1$}}}
\end{picture}

%% file: nonTij.pstex_t
\begin{picture}(0,0)%
\includegraphics{nonTij.pstex}%
\end{picture}%
\setlength{\unitlength}{3947sp}%
\begingroup\makeatletter\ifx\SetFigFont\undefined
\def\x#1#2#3#4#5#6#7\relax{\def\x{#1#2#3#4#5#6}}%
\expandafter\x\fmtname xxxxxx\relax \def\y{splain}%
\ifx\x\y   
\gdef\SetFigFont#1#2#3{%
  \ifnum #1<17\tiny\else \ifnum #1<20\small\else
  \ifnum #1<24\normalsize\else \ifnum #1<29\large\else
  \ifnum #1<34\Large\else \ifnum #1<41\LARGE\else
     \huge\fi\fi\fi\fi\fi\fi
  \csname #3\endcsname}%
\else
\gdef\SetFigFont#1#2#3{\begingroup
  \count@#1\relax \ifnum 25<\count@\count@25\fi
  \def\x{\endgroup\@setsize\SetFigFont{#2pt}}%
  \expandafter\x
    \csname \romannumeral\the\count@ pt\expandafter\endcsname
    \csname @\romannumeral\the\count@ pt\endcsname
  \csname #3\endcsname}%
\fi
\fi\endgroup
\begin{picture}(2185,2748)(663,-2437)
\put(946,-1006){\makebox(0,0)[lb]{\smash{\SetFigFont{12}{14.4}{rm}$P_{1}$}}}
\put(2566,-1006){\makebox(0,0)[lb]{\smash{\SetFigFont{12}{14.4}{rm}$P_{n}$}}}
\put(1531,-1006){\makebox(0,0)[lb]{\smash{\SetFigFont{12}{14.4}{rm}$P_{i}$}}}
\put(1711,-2356){\makebox(0,0)[lb]{\smash{\SetFigFont{17}{20.4}{rm}$T_{i,j}$}}}
\put(2116,-961){\makebox(0,0)[lb]{\smash{\SetFigFont{12}{14.4}{rm}$P_{j}$}}}
\end{picture}

%% file: nonAjr.pstex_t
\begin{picture}(0,0)%
\includegraphics{nonAjr.pstex}%
\end{picture}%
\setlength{\unitlength}{3947sp}%
\begingroup\makeatletter\ifx\SetFigFont\undefined
\def\x#1#2#3#4#5#6#7\relax{\def\x{#1#2#3#4#5#6}}%
\expandafter\x\fmtname xxxxxx\relax \def\y{splain}%
\ifx\x\y   
\gdef\SetFigFont#1#2#3{%
  \ifnum #1<17\tiny\else \ifnum #1<20\small\else
  \ifnum #1<24\normalsize\else \ifnum #1<29\large\else
  \ifnum #1<34\Large\else \ifnum #1<41\LARGE\else
     \huge\fi\fi\fi\fi\fi\fi
  \csname #3\endcsname}%
\else
\gdef\SetFigFont#1#2#3{\begingroup
  \count@#1\relax \ifnum 25<\count@\count@25\fi
  \def\x{\endgroup\@setsize\SetFigFont{#2pt}}%
  \expandafter\x
    \csname \romannumeral\the\count@ pt\expandafter\endcsname
    \csname @\romannumeral\the\count@ pt\endcsname
  \csname #3\endcsname}%
\fi
\fi\endgroup
\begin{picture}(2655,2682)(361,-2098)
\put(2161,389){\makebox(0,0)[lb]{\smash{\SetFigFont{12}{14.4}{rm}$s_{i,r}$}}}
\put(2476,-1771){\makebox(0,0)[lb]{\smash{\SetFigFont{12}{14.4}{rm}$e_i$}}}
\put(406,-1276){\makebox(0,0)[lb]{\smash{\SetFigFont{12}{14.4}{rm}$P_i$}}}
\put(1216,389){\makebox(0,0)[lb]{\smash{\SetFigFont{12}{14.4}{rm}$s_{i,r}$}}}
\put(2161,-1096){\makebox(0,0)[lb]{\smash{\SetFigFont{12}{14.4}{rm}$P_n$}}}
\put(1756,-1456){\makebox(0,0)[lb]{\smash{\SetFigFont{12}{14.4}{rm}$P_j$}}}
\put(2476,-1096){\makebox(0,0)[lb]{\smash{\SetFigFont{12}{14.4}{rm}$P_1$}}}
\put(361,-871){\makebox(0,0)[lb]{\smash{\SetFigFont{12}{14.4}{rm}$s_{i,1}$}}}
\put(3016,-826){\makebox(0,0)[lb]{\smash{\SetFigFont{12}{14.4}{rm}$s_{i,g}$}}}
\put(1036,-1771){\makebox(0,0)[lb]{\smash{\SetFigFont{12}{14.4}{rm}$e_i$}}}
\end{picture}